\documentclass[final]{siamltex}%
\usepackage{algpseudocode}
\usepackage{amsmath}
\usepackage{amsfonts}
\usepackage{amssymb}
\usepackage[draft]{graphicx}
\usepackage[mathscr]{euscript}
\usepackage{multirow}
\usepackage{subfigure}
\usepackage{url}
\usepackage[pdftex]{thumbpdf}
\usepackage{verbatim}
\usepackage[pdftex,
pdfstartview=FitH,bookmarks=true,bookmarksnumbered=true,bookmarksopen=true,hypertexnames=false,breaklinks=true,
colorlinks=true,linkcolor=blue,anchorcolor=blue,
citecolor=blue,filecolor=blue,menucolor=blue]{hyperref}%
\providecommand{\U}[1]{\protect\rule{.1in}{.1in}}
\newtheorem{algorithm}[theorem]{Algorithm}

\newtheorem{remark}[theorem]{Remark}
\begin{document}

\title{Stochastic Galerkin Methods for the Steady-State Navier-Stokes Equations
\thanks{This work is based upon work supported by the U.\, S.\, Department of
Energy Office of Advanced Scientific Computing Research, Applied Mathematics
program under Award Number DE-SC0009301, and by the U.\, S.\, National Science
Foundation under grants DMS1418754 and DMS1521563.}}
\author{Bed\v{r}ich Soused\'{\i}k\thanks{Department of Mathematics and Statistics,
University of Maryland, Baltimore County, 1000 Hilltop Circle, Baltimore,
MD~21250 (\texttt{sousedik@umbc.edu}). }
\and Howard C. Elman\thanks{Department of Computer Science and Institute for
Advanced Computer Studies, University of Maryland, College Park, MD 20742
(\texttt{elman@cs.umd.edu}) }}
\maketitle

\begin{abstract}
We study the steady-state Navier-Stokes equations in the context of stochastic
finite element discretizations. Specifically, we assume that the viscosity is
a random field given in the form of a generalized polynomial chaos expansion.
For the resulting stochastic problem, we formulate the model and linearization
schemes using Picard and Newton iterations in the framework of the stochastic
Galerkin method, and we explore properties of the resulting stochastic
solutions. We also propose a preconditioner for solving the linear systems of
equations arising at each step of the stochastic (Galerkin) nonlinear
iteration and demonstrate its effectiveness for solving a set of benchmark problems.

\end{abstract}

\section{Introduction}

\label{sec:introduction} Models of mathematical physics are typically based on
partial differential equations (PDEs) that use parameters as input data. In
many situations, the values of parameters are not known precisely and are
modeled as random fields, giving rise to \emph{stochastic partial differential
equations}. In this study we focus on models from fluid dynamics, in
particular the stochastic Stokes and the Navier-Stokes equations. We consider
the viscosity as a random field modeled as colored noise, and we use numerical
methods based on spectral methods, specifically, the generalized polynomial
chaos (gPC)
framework~\cite{Ghanem-1991-SFE,LeMaitre-2010-SMU,Xiu-2010-NMS,Xiu-2002-WAP}.
That is, the viscosity is given by a gPC expansion, and we seek gPC expansions
of the velocity and pressure solutions.

There is a number of reasons to motivate our interest in Navier-Stokes
equations with stochastic viscosity. For example, the exact value of viscosity
may not be known, due to measurement error, the presence of contaminants with
uncertain concentrations, or of multiple phases with uncertain ratios.
Alternatively, the fluid properties might be influenced by an external field,
with applications for example in magnetohydrodynamics. Specifically, we assume
that the viscosity$~\nu$ depends on a set of random variables $\xi$. This
means that the Reynolds number,
\[
\operatorname{Re}\left(  \xi\right)  =\frac{UL}{\nu\left(  \xi\right)  },
\]
where $U$ is the characteristic velocity and $L$ is the characteristic length,
is also stochastic. Consequently, the solution variables are random fields,
and different realizations of the viscosity give rise to realizations of the
velocities and pressures. As observed in~\cite{Powell-2012-PSS}, there are
other possible formulations and interpretations of fluid flow with stochastic
Reynolds number for example, where the velocity is fixed but the volume of
fluid moving into a channel is uncertain so the uncertainty derives from the
Dirichlet inflow boundary condition.

We consider models of steady-state stochastic motion of an incompressible
fluid moving in a domain~$D\subset\mathbb{R}^{2}$. Extension to
three-dimensional models is straightforward. We formulate the stochastic
Stokes and Navier-Stokes equations using the stochastic finite element method,
assuming that the viscosity has a general probability distribution
parametrized by a gPC\ expansion. We describe linearization schemes based on
Picard and Newton iteration for the stochastic Galerkin method, and we explore
properties of the solutions obtained, including a comparison of the stochastic
Galerkin solutions with those obtained using other approaches, such as Monte
Carlo and stochastic collocation methods~\cite{Xiu-2010-NMS}. Finally, we
propose efficient hierarchical preconditioners for iterative solution of the
linear systems solved at each step of the nonlinear iteration in the context
of the stochastic Galerkin method. Our approach is related to recent work by
Powell and Silvester~\cite{Powell-2012-PSS}. However, besides using a general
parametrization of the viscosity, our formulation of the stochastic Galerkin
system allows straightforward application of state-of-the-art deterministic
preconditioners by wrapping them in the hierarchical preconditioner developed
in~\cite{Sousedik-2014-THP}. For alternative preconditioners see,
e.g.,~\cite{Brezina-2014-SAA,Giraldi-2014-LAI,Pellissetti-2000-ISS,
Powell-2009-BDP,Rosseel-2010-ISS,Sousedik-2014-HSC,Ullmann-2010-KPD}. Finally,
we note that there exist related approaches based on stochastic perturbation
methods~\cite{Kaminski-2013-SPM}, important developments also include
reduced-order models such as~\cite{Elman-2013-RBC,Tamellini-2014-MRB}, and an
overview of existing methods for stochastic computational fluid dynamics can
be found in the monograph~\cite{LeMaitre-2010-SMU}.

The paper is organized as follows. In Section~\ref{sec:NS}, we recall the
deterministic steady-state Navier-Stokes equations and their discrete form. In
Section~\ref{sec:stochastic}, we formulate the model with stochastic
viscosity, derive linearization schemes for the stochastic Galerkin
formulation of the model, and explore properties of the resulting solutions
for a set of benchmark problems that model the flow over an obstacle. In
Section~\ref{sec:preconditioner} we introduce a preconditioner for the
stochastic Galerkin linear systems solved at each step of the nonlinear
iteration, and in Section~\ref{sec:conclusion} we summarize our work.

\section{Deterministic Navier-Stokes equations}

\label{sec:NS} We begin by defining the model and notation,
following~\cite{Elman-2014-FEF}. For the deterministic Navier-Stokes
equations, we wish to find velocity$~\vec{u}$ and pressure$~p$ such that
\begin{align}
-\nu\nabla^{2}\vec{u}+\left(  \vec{u}\cdot\nabla\right) \vec{u}+\nabla p &
=\vec{f},\label{eq:NS-1}\\
\nabla\cdot\vec{u} &  =0,\label{eq:NS-2}%
\end{align}
in a spatial domain$~D$, satisfying boundary conditions%
\begin{align}
\vec{u} &  =\vec{g},\quad\text{on }\Gamma_{\text{Dir}},\label{eq:Dir-bc}\\
\nu\nabla\vec{u}\cdot\vec{n}-p\vec{n} &  =\vec{0},\quad\text{on }%
\Gamma_{\text{Neu}},\label{eq:Neu-bc}%
\end{align}
where$~\partial D=\overline{\Gamma}_{\text{Dir}}\cup\overline{\Gamma
}_{\text{Neu}}$, and assuming sufficient regularity of the data. Dropping the
convective term $\left(  \vec{u}\cdot\nabla\right)  \vec{u}$
from~(\ref{eq:NS-1}) yields the Stokes problem
\begin{align}
-\nu\nabla^{2}\vec{u}+\nabla p &  =\vec{f},\label{eq:S-1}\\
\nabla\cdot\vec{u} &  =0.\label{eq:S-2}%
\end{align}
The mixed variational formulation of~(\ref{eq:NS-1})--(\ref{eq:NS-2}) is to
find $\left(  \vec{u},p\right)  \in\left(  V_{E},Q_{D}\right)  $ such that
\begin{align}
\nu\int_{D}\nabla\vec{u}:\nabla\vec{v}+\int_{D}\left(  \vec{u}\cdot\nabla
\vec{u}\right) \vec{v}-\int_{D}p\left(  \nabla\cdot\vec{v}\right)   &
=\int_{D}\vec{f}\cdot\vec{v},\quad\forall\vec{v}\in V_{D}%
,\label{eq:NS-variational-1}\\
\int_{D}q\left(  \nabla\cdot\vec{u}\right)   &  =0,\quad\forall q\in
Q_{D},\label{eq:NS-variational-2}%
\end{align}
where $\left(  V_{D},Q_{D}\right)  $ is a pair of spaces satisfying the
inf-sup condition and $V_{E}$ is an extension of $V_{D}$ containing velocity
vectors that satisfy the Dirichlet boundary
conditions~\cite{Brezzi-1991-MHF,Elman-2014-FEF,Girault-1986-FEM}.

Let $c(\vec{z};\vec{u},\vec{v})\equiv\int_{\Omega}\left(  \vec{z}\cdot
\nabla\vec{u}\right)  \cdot\vec{v}$. Because the
problem~(\ref{eq:NS-variational-1})--(\ref{eq:NS-variational-2}) is nonlinear,
it is solved using a linearization scheme in the form of Newton or Picard
iteration, derived as follows. Consider the solution $\left(  \vec
{u},p\right)  $ of~(\ref{eq:NS-variational-1})--(\ref{eq:NS-variational-2}) to
be given as $\vec{u}=\vec{u}^{n}+\delta\vec{u}^{n}$ and $p=p^{n}+\delta p^{n}%
$. Substituting into~(\ref{eq:NS-variational-1})--(\ref{eq:NS-variational-2})
and neglecting the quadratic term $c(\delta\vec{u}^{n};\delta\vec{u}^{n}%
,\vec{v})$ gives
\begin{align}
\nu\int_{D}\nabla\delta\vec{u}^{n}:\nabla\vec{v}+c(\delta\vec{u}^{n};\vec
{u}^{n},\vec{v})+c(\vec{u}^{n};\delta\vec{u}^{n},\vec{v})-\int_{D}\delta
p^{n}\left(  \nabla\cdot\vec{v}\right)   &  =R^{n}\left(  \vec{v}\right)
,\label{eq:Newton-1}\\
\int_{D}q\left(  \nabla\cdot\delta\vec{u}^{n}\right)   &  =r^{n}\left(
q\right)  ,\label{eq:Newton-2}%
\end{align}
where
\begin{align}
R^{n}\left(  \vec{v}\right)   &  =\int_{D}\vec{f}\cdot\vec{v}-\nu\int%
_{D}\nabla\vec{u}^{n}:\nabla\vec{v}-c(\vec{u}^{n};\vec{u}^{n},\vec{v}%
)+\int_{D}p^{n}\left(  \nabla\cdot\vec{v}\right)  ,\label{eq:residual-1}\\
r^{n}\left(  q\right)   &  =-\int_{D}q\left(  \nabla\cdot\vec{u}^{n}\right)
.\label{eq:residual-2}%
\end{align}
Step~$n$\ of the \emph{Newton iteration} obtains $\left(  \delta\vec{u}%
^{n},\delta p^{n}\right)  $ from~(\ref{eq:Newton-1})--(\ref{eq:Newton-2}) and
updates the solution as
\begin{align}
\vec{u}^{n+1} &  =\vec{u}^{n}+\delta\vec{u}^{n},\label{eq:Newton-increment-1}%
\\
p^{n+1} &  =p^{n}+\delta p^{n}.\label{eq:Newton-increment-2}%
\end{align}
Step~$n$\ of the \emph{Picard iteration} omits the term $c(\delta\vec{u}%
^{n};\vec{u}^{n},\vec{v})$ in (\ref{eq:Newton-1}), giving
\begin{align}
\nu\int_{D}\nabla\delta\vec{u}^{n}:\nabla\vec{v}+c(\vec{u}^{n};\delta\vec
{u}^{n},\vec{v})-\int_{D}\delta p^{n}\left(  \nabla\cdot\vec{v}\right)   &
=R^{n}\left(  \vec{v}\right)  ,\label{eq:Picard-1}\\
\int_{D}q\left(  \nabla\cdot\delta\vec{u}^{n}\right)   &  =r^{n}\left(
q\right)  .\label{eq:Picard-2}%
\end{align}

Consider the discretization of~(\ref{eq:NS-1})--(\ref{eq:NS-2}) by a
div-stable mixed finite element method; for experiments discussed below, we
used Taylor-Hood elements~\cite{Elman-2014-FEF}. Let the bases for the
velocity and pressure spaces be denoted $\left\{  \phi_{i}\right\}
_{i=1}^{N_{u}}$\ and $\left\{  \varphi_{i}\right\}  _{i=1}^{N_{p}}$,
respectively. In matrix terminology, each nonlinear iteration entails solving
a linear system
\begin{equation}
\left[
\begin{array}
[c]{cc}%
\mathbf{F}^{n} & \mathbf{B}^{T}\\
\mathbf{B} & \mathbf{0}%
\end{array}
\right]  \left[
\begin{array}
[c]{c}%
\delta\mathbf{u}^{n}\\
\delta\mathbf{p}^{n}%
\end{array}
\right]  =\left[
\begin{array}
[c]{c}%
\mathbf{R}^{n}\\
\mathbf{r}^{n}%
\end{array}
\right]  ,\label{eq:Newton}%
\end{equation}
followed by an update of the solution
\begin{align}
\mathbf{u}^{n+1} &  =\mathbf{u}^{n}+\delta\mathbf{u}^{n},\label{eq:update-1}\\
\mathbf{p}^{n+1} &  =\mathbf{p}^{n}+\delta\mathbf{p}^{n}.\label{eq:update-2}%
\end{align}
For Newton's method, $\mathbf{F}^{n}$ is the Jacobian matrix, a sum of the
vector-Laplacian matrix~$\mathbf{A}$, the vector-convection matrix~$\mathbf{N}%
^{n}$, and the Newton derivative matrix~$\mathbf{W}^{n}$,
\begin{equation}
\mathbf{F}^{n}=\mathbf{A}+\mathbf{N}^{n}+\mathbf{W}^{n}%
,\label{Newton-Jacobian}%
\end{equation}
where
\begin{align*}
\mathbf{A}\mathbf{=}\left[  a_{ab}\right]  ,\qquad &  a_{ab}=\nu\int%
_{D}\,\nabla\phi_{b}:\nabla\phi_{a}\mathbf{,}\\
\mathbf{N}^{n}=\left[  n_{ab}^{n}\right]  ,\qquad &  n_{ab}^{n}=\int%
_{D}\left(  u^{n}\cdot\nabla\phi_{b}\right)  \cdot\phi_{a},\\
\mathbf{W}^{n}=\left[  w_{ab}^{n}\right]  ,\qquad &  w_{ab}^{n}=\int%
_{D}\left(  \phi_{b}\cdot\nabla u^{n}\right)  \cdot\phi_{a}.
\end{align*}
For Picard iteration, the Newton derivative matrix~$\mathbf{W}^{n}$ is
dropped, and $\mathbf{F}^{n}=\mathbf{A}+\mathbf{N}^{n}$. The divergence
matrix~$B$ is defined as
\begin{equation}
\mathbf{B}=\left[  b_{cd}\right]  ,\qquad b_{cd}=\int_{D}\phi_{d}\left(
\nabla\cdot\varphi_{c}\right)  .\label{eq:B}%
\end{equation}
The residuals at step~$n$ of both nonlinear iterations are computed as
\begin{equation}
\left[
\begin{array}
[c]{c}%
\mathbf{R}^{n}\\
\mathbf{r}^{n}%
\end{array}
\right]  =\left[
\begin{array}
[c]{c}%
\mathbf{f}\\
\mathbf{g}%
\end{array}
\right]  -\left[
\begin{array}
[c]{cc}%
\mathbf{P}^{n} & \mathbf{B}^{T}\\
\mathbf{B} & \mathbf{0}%
\end{array}
\right]  \left[
\begin{array}
[c]{c}%
\mathbf{u}^{n}\\
\mathbf{p}^{n}%
\end{array}
\right]  ,\label{eq:residual}%
\end{equation}
where $\mathbf{P}^{n}=\mathbf{A}+\mathbf{N}^{n}$ and $\mathbf{f}$ is a
discrete version of the forcing function of~(\ref{eq:NS-1}%
).\footnote{Throughout this study, we use the convention that the right-hand
sides of discrete systems incorporate Dirichlet boundary data for velocities.}

\section{The Navier-Stokes equations with stochastic viscosity}

\label{sec:stochastic} Let $\left(  \Omega,\mathcal{F},\mathcal{P}\right)  $
represent a complete probability space, where $\Omega$ is the sample space,
$\mathcal{F}$ is a$~\sigma$-algebra on~$\Omega$ and $\mathcal{P}$ is a
probability measure. We will assume that the randomness in the model is
induced by a vector of independent, identically distributed (i.i.d.)\ random
variables $\xi=\left(  \xi_{1},\dots,\xi_{N}\right)  ^{T}$ such that
$\xi:\Omega\rightarrow\Gamma\subset\mathbb{R}^{N}$. Let $\mathcal{F}_{\xi
}\mathcal{\subset F}$\ denote the $\sigma$-algebra generated by$~\xi$, and let
$\mu\left(  \xi\right)  $ denote the joint probability density measure
for~$\xi$. The expected value of the product of random variables$~u$
and$~v$\ that depend on$~\xi$\ determines a Hilbert space~$T_{\Gamma}\equiv
L^{2}\left(  \Omega,\mathcal{F}_{\xi},\mu\right)  $ with inner product
\begin{equation}
\left\langle u,v\right\rangle =\mathbb{E}\left[  uv\right]  =\int_{\Gamma
}u\left(  \xi\right)  v\left(  \xi\right)  \,d\mu\left(  \xi\right)  ,
\label{eq:E}%
\end{equation}
where the symbol~$\mathbb{E}$ denotes mathematical expectation.

\subsection{The stochastic Galerkin formulation}

\label{sec:Galerkin} The counterpart of the variational
formulation~(\ref{eq:NS-variational-1})--(\ref{eq:NS-variational-2}) consists
of performing a Galerkin projection on the space $T_{\Gamma}$ using
mathematical expectation in the sense of~(\ref{eq:E}). That is, we seek the
velocity~$\vec{u}$, a random field in~$V_{E}\otimes T_{\Gamma}$, and the
pressure $p\in Q_{D}\otimes T_{\Gamma}$, such that
\begin{align}
\mathbb{E}\left[  \int_{D}\!\nu\,\nabla\vec{u}:\nabla\vec{v}+\int_{D}\!\left(
\vec{u}\cdot\nabla\vec{u}\right)  \vec{v}-\int_{D}\!p\left(  \nabla\cdot
\vec{v}\right)  \right]   &  =\mathbb{E}\left[  \int_{D}\!\vec{f}\cdot\vec
{v}\right]  \quad\forall\vec{v}\in V_{D}\otimes T_{\Gamma},\label{eq:NS-1-s}\\
\mathbb{E}\left[  \int_{D}\!q\left(  \nabla\cdot\vec{u}\right)  \right]   &
=0\quad\forall q\in Q_{D}\otimes T_{\Gamma}.\label{eq:NS-2-s}%
\end{align}
The stochastic counterpart of the Newton iteration~(\ref{eq:Newton-1}%
)--(\ref{eq:Newton-2}) is
\begin{align}
\mathbb{E}\left[  \int_{D}\nu\,\nabla\delta\vec{u}^{n}:\nabla\vec{v}+c(\vec
{u}^{n};\delta\vec{u}^{n},\vec{v})+c(\delta\vec{u}^{n};\vec{u}^{n},\vec
{v})-\int_{D}\delta p^{n}\left(  \nabla\cdot\vec{v}\right)  \right]   &
=R^{n},\label{eq:Newton-1-s}\\
\mathbb{E}\left[  \int_{D}q\left(  \nabla\cdot\delta\vec{u}^{n}\right)
\right]   &  =r^{n},\label{eq:Newton-2-s}%
\end{align}
where
\begin{align}
R^{n}\left(  \vec{v}\right)   &  =\mathbb{E}\left[  \int_{D}\vec{f}\cdot
\vec{v}-\int_{D}\nu\,\nabla\vec{u}^{n}:\nabla\vec{v}-c(\vec{u}^{n};\vec{u}%
^{n},\vec{v})+\int_{D}p^{n}\left(  \nabla\cdot\vec{v}\right)  \right]
,\label{eq:s-residual-1}\\
r^{n}\left(  q\right)   &  =-\mathbb{E}\left[  \int_{D}q\left(  \nabla
\cdot\vec{u}^{n}\right)  \right]  .\label{eq:s-residual-2}%
\end{align}
The analogue for Picard iteration omits $c(\delta\vec{u}^{n};\vec{u}^{n}%
,\vec{v})$ from~(\ref{eq:Newton-1-s}):
\begin{equation}
\mathbb{E}\left[  \int_{D}\nu\,\nabla\delta\vec{u}^{n}:\nabla\vec{v}+c(\vec
{u}^{n};\delta\vec{u}^{n},\vec{v})-\int_{D}\delta p^{n}\left(  \nabla\cdot
\vec{v}\right)  \right]  =R^{n}.\label{eq:Picard-s}%
\end{equation}

In computations, we will use a finite-dimensional subspace $T_{P}\subset
T_{\Gamma}$ spanned by a set of polynomials $\left\{  \psi_{\ell}\left(
\xi\right)  \right\}  $ that are orthogonal with respect to the density
function~$\mu$, that is $\left\langle \psi_{k},\psi_{\ell}\right\rangle
=\delta_{k\ell}$. This is referred to as the gPC basis;
see~\cite{Ghanem-1991-SFE,Xiu-2002-WAP} for details and discussion.
For~$T_{P}$, we will use the space spanned by multivariate polynomials
in~$\{\xi_{j}\}_{j=1}^{N}$ of total degree~$P$, which has
dimension~$M={\renewcommand{\arraystretch}{0.8}\left(  \!\!%
\begin{array}
[c]{c}%
N+P\\
P
\end{array}
\!\!\right)  }$. We will also assume that the viscosity is given by a gPC
expansion
\begin{equation}
\nu=\sum_{\ell=0}^{M_{\nu}-1}\nu_{\ell}\left(  x\right)  \psi_{\ell}\left(
\xi\right)  ,\label{eq:viscosity-gPC}%
\end{equation}
where $\left\{  \nu_{\ell}\left(  x\right)  \right\}  $ is a set of given
deterministic spatial functions.

\subsection{Stochastic Galerkin finite element formulation}

\label{sec:Galerkin-FEM} We discretize~(\ref{eq:Newton-1-s})
(or~(\ref{eq:Picard-s})) and~(\ref{eq:Newton-2-s}) using div-stable finite
elements as in Section~\ref{sec:NS} together with the gPC basis for~$T_{P}$.
For simplicity, we assume that the right-hand side~$f\left(  x\right)  $ and
the Dirichlet boundary conditions~(\ref{eq:Dir-bc}) are deterministic. This
means in particular that, as in the deterministic case, the boundary
conditions can be incorporated into right-hand side vectors (specified
as~$\mathbf{y}$ below). Thus, we seek a discrete approximation of the solution
of the form
\begin{align}
\vec{u}\left(  x,\xi\right)   &  \approx\sum_{k=0}^{M-1}\sum_{i=1}^{N_{u}%
}u_{ik}\phi_{i}(x)\psi_{k}(\xi)=\sum_{k=0}^{M-1}\vec{u}_{k}(x)\psi_{k}%
(\xi),\label{eq:u-gPC}\\
p\left(  x,\xi\right)   &  \approx\sum_{k=0}^{M-1}\sum_{j=1}^{N_{p}}%
p_{jk}\varphi_{j}(x)\psi_{k}(\xi)=\sum_{k=0}^{M-1}p_{k}(x)\psi_{k}%
(\xi),\label{eq:p-gPC}%
\end{align}

The structure of the discrete operators depends on the ordering of the unknown
coefficients~$\{u_{ik}\}$, $\{p_{jk}\}$. We will group velocity-pressure pairs
for each~$k$, the index of stochastic basis functions (and order equations in
the same way), giving the ordered list of coefficients
\begin{equation}
u_{1:N_{u},0},p_{1:N_{p},0},\,u_{1:N_{u},1},p_{1:N_{p},1},\,\ldots
,\,u_{1:N_{u},M-1},p_{1:N_{p},M-1}.\label{eq:coefficient-order}%
\end{equation}
To describe the discrete structure, we first consider the stochastic version
of the Stokes problem~(\ref{eq:S-1})--(\ref{eq:S-2}), where the convection
term~$c(\cdot;\cdot,\cdot)$ is not present in~(\ref{eq:Newton-1-s})
and~(\ref{eq:Picard-s}). The discrete stochastic Stokes operator is built from
the discrete components of the vector-Laplacian
\begin{equation}
\mathbf{A}_{\ell}\mathbf{=}\left[  a_{\ell,ab}\right]  ,\quad a_{\ell
,ab}=\left(  \int_{D}\nu_{\ell}\left(  x\right)  \,\nabla\phi_{b}:\nabla
\phi_{a}\right)  ,\qquad\ell=1,\ldots,M_{\nu}%
-1,\label{stochastic-vector-Laplacian}%
\end{equation}
which are incorporated into the block matrices
\begin{equation}
\mathcal{S}_{0}=\left[
\begin{array}
[c]{cc}%
\mathbf{A}_{0} & \mathbf{B}^{T}\\
\mathbf{B} & \mathbf{0}%
\end{array}
\right]  ,\qquad\mathcal{S}_{\ell}=\left[
\begin{array}
[c]{cc}%
\mathbf{A}_{\ell} & \mathbf{0}\\
\mathbf{0} & \mathbf{0}%
\end{array}
\right]  ,\quad\ell=1,\dots,M_{\nu}-1.\label{stochastic-saddle-point}%
\end{equation}
These operators will be coupled with matrices arising from terms in~$T_{P}$,
\begin{equation}
\mathbf{H}_{\ell}=\left[  h_{\ell,jk}\right]  ,\quad h_{\ell,jk}%
\equiv\mathbb{E}\left[  \psi_{\ell}\psi_{j}\psi_{k}\right]  ,\qquad
\ell=0,\dots,M_{\nu}-1,\quad j,k=0,\dots
,M-1.\label{eq:stochastic-matrix-forms}%
\end{equation}
Combining the expressions from~(\ref{stochastic-vector-Laplacian}),
(\ref{stochastic-saddle-point}) and~(\ref{eq:stochastic-matrix-forms}) and
using the ordering~(\ref{eq:coefficient-order}) gives the discrete stochastic
Stokes system
\begin{equation}
\left(  \sum_{\ell=0}^{M_{\nu}-1}\mathbf{H}_{\ell}\otimes\mathcal{S}_{\ell
}\right)  \mathbf{v}=\mathbf{y},\label{eq:Stokes-s}%
\end{equation}
where $\otimes$ corresponds to the matrix Kronecker product. The unknown
vector $\mathbf{v}$ corresponds to the ordered list of coefficients
in~(\ref{eq:coefficient-order}) and the right-hand side is ordered in an
analogous way. Note that~$\mathbf{H}_{0}$ is the identity matrix of order~$M$.

\begin{remark}
\label{rem:order} With this ordering, the coefficient matrix contains a set
of~$M$ block~$2\times2$ matrices of saddle-point structure along its block
diagonal, given by
\[
\mathcal{S}_{0}+\sum_{\ell=1}^{M_{\nu}-1}h_{\ell,jj}\mathcal{S}_{\ell},\qquad
j=0,\ldots,M-1.
\]
This enables the use of existing deterministic solvers for the individual
diagonal blocks. An alternative ordering based on the blocking of \emph{all}
velocity coefficients followed by all pressure coefficients, considered
in~\cite{Powell-2012-PSS}, produces a matrix of global saddle-point structure.
\end{remark}

The matrices arising from the linearized stochastic Navier-Stokes equations
augment the Stokes systems with stochastic variants of the vector-convection
matrix and Newton derivative matrix appearing in~(\ref{Newton-Jacobian}). In
particular, at step~$n$ of the nonlinear iteration, let~$\vec{u}_{\ell}%
^{n}(x)$ be the $\ell$th term of the velocity iterate (as in the expression on
the right in~(\ref{eq:u-gPC}) for $k=\ell$), and let
\begin{align*}
\mathbf{N}_{\ell}^{n}=\left[  n_{\ell,ab}^{n}\right]  ,\qquad &  n_{\ell
,ab}^{n}=\int_{D}\left(  \vec{u}_{\ell}^{n}\cdot\nabla\phi_{b}\right)
\cdot\phi_{a},\\
\mathbf{W}_{\ell}^{n}=\left[  w_{\ell,ab}^{n}\right]  ,\qquad &  w_{\ell
,ab}^{n}=\int_{D}\left(  \phi_{b}\cdot\nabla\vec{u}_{\ell}^{n}\right)
\cdot\phi_{a}.
\end{align*}
Then the analogues of~(\ref{stochastic-vector-Laplacian}%
)--(\ref{stochastic-saddle-point}) are
\begin{align}
\mathbf{F}_{\ell}^{n}  &  =\mathbf{A}_{\ell}+\mathbf{N}_{\ell}^{n}%
+\mathbf{W_{\ell}}^{n}, &  &
\mbox{for  the stochastic Newton method}\label{eq:stoch-Newton}\\
\mathbf{F}_{\ell}^{n}  &  =\mathbf{A}_{\ell}+\mathbf{N}_{\ell}^{n}, &  &
\mbox{for  stochastic Picard iteration,\ \,} \label{eq:stoch-Picard}%
\end{align}
so for Newton's method
\begin{equation}
\mathcal{F}_{0}^{n}=\left[
\begin{array}
[c]{cc}%
\mathbf{F}_{0}^{n} & \mathbf{B}^{T}\\
\mathbf{B} & \mathbf{0}%
\end{array}
\right]  ,\qquad\mathcal{F}_{\ell}^{n}=\left[
\begin{array}
[c]{cc}%
\mathbf{F}_{\ell} & \mathbf{0}\\
\mathbf{0} & \mathbf{0}%
\end{array}
\right]  , \label{eq:linearized-operator}%
\end{equation}
and as above, for Picard iteration the Newton derivative
matrices~$\{\mathbf{W}^{n}\}$ are dropped. Note that $\ell=0,\ldots
,\widehat{M}-1$ here, where $\widehat{M}=\max\left(  M,M_{\nu}\right)  $. (In
particular, if $M_{\nu}>M$, we set $\mathbf{N}_{\ell}^{n}=\mathbf{W}_{\ell
}^{n}=0$ for $\ell=M+1,\dots,M_{\nu}-1$.) Step~$n$ of the stochastic nonlinear
iteration entails solving a linear system and updating,
\begin{equation}
\left[  \sum_{\ell=0}^{\widehat{M}-1}\mathbf{H}_{\ell}\otimes\mathcal{F}%
_{\ell}^{n}\right]  \delta\mathbf{v}^{n}=\mathcal{R}^{n},\qquad\mathbf{v}%
^{n+1}=\mathbf{v}^{n}+\delta\mathbf{v}^{n}, \label{eq:linearized-system}%
\end{equation}
where
\begin{equation}
\mathcal{R}^{n}=\mathbf{y}-\left[  \sum_{\ell=0}^{\widehat{M}-1}%
\mathbf{H}_{\ell}\otimes\mathcal{P}_{\ell}^{n}\right]  \mathbf{v}^{n},
\label{eq:R_gPC}%
\end{equation}
$\mathbf{v}^{n}$ and $\delta\mathbf{v}^{n}$ are vectors of current velocity
and pressure coefficients and updates, respectively, ordered as
in~(\ref{eq:coefficient-order}), $\mathbf{y}$ is the similarly ordered
right-hand side determined from the forcing function and Dirichlet boundary
data, and
\[
\mathcal{P}_{0}^{n}=\left[
\begin{array}
[c]{cc}%
\mathbf{A}_{0}+\mathbf{N}_{0}^{n} & \mathbf{B}^{T}\\
\mathbf{B} & \mathbf{0}%
\end{array}
\right]  ,\qquad\mathcal{P}_{\ell}^{n}=\left[
\begin{array}
[c]{cc}%
\mathbf{A}_{\ell}+\mathbf{N}_{\ell}^{n} & \mathbf{0}\\
\mathbf{0} & \mathbf{0}%
\end{array}
\right]  ;
\]
note that the~$(1,1)$-blocks here are as in~(\ref{eq:stoch-Picard}).

\subsection{Sampling methods}

\label{sec:collocation} In experiments described below, we compare some
results obtained using stochastic Galerkin methods to those obtained from
Monte Carlo and stochastic collocation. We briefly describe these approaches here.

Both Monte Carlo and stochastic collocation methods are based on sampling.
This entails the solution of a number of mutually independent deterministic
problems at a set of sample points~$\left\{  \xi^{\left(  q\right)  }\right\}
$, which give realizations of the viscosity~(\ref{eq:viscosity-gPC}). That is,
a realization of viscosity~$\nu\left(  \xi^{\left(  q\right)  }\right)  $
gives rise to deterministic functions~$\vec{u}\left(  \cdot,\xi^{\left(
q\right)  }\right)  $ and$~p\left(  \cdot,\xi^{\left(  q\right)  }\right)  $
on~$D$ that satisfy the standard deterministic Navier-Stokes equations, and to
finite-element approximations~$\vec{u}^{(q)}(x)$,~$p^{(q)}(x)$.

In the Monte Carlo method, the$~N_{MC}$ sample points are generated randomly,
following the distribution of the random variables$~\xi$, and moments of the
solution are obtained from ensemble averaging.
For stochastic collocation, the sample points consist of a set of
predetermined \emph{collocation points}. This approach derives from a
methodology for performing quadrature or interpolation in multidimensional
space using a small number of points, a so-called sparse
grid~\cite{Gerstner-1998-NIU,Novak-1996-HDI}. There are two ways to implement
stochastic collocation to obtain the coefficients in~(\ref{eq:u-gPC}%
)--(\ref{eq:p-gPC}), either by constructing a Lagrange interpolating
polynomial, or, in the so-called pseudospectral approach, by performing a
discrete projection into~$T_{P}$~\cite{Xiu-2010-NMS}. We use the second
approach because it facilitates a direct comparison with the stochastic
Galerkin method. In particular, the coefficients
are determined using a quadrature
\[
u_{ik}=\sum_{q=1}^{N_{q}}\vec{u}^{(q)}\left(  x_{i}\right)  \,\psi_{k}\left(
\xi^{\left(  q\right)  }\right)  \,w^{\left(  q\right)  },\qquad p_{ik}%
=\sum_{q=1}^{N_{q}}p^{(q)}\left(  x_{i}\right)  \,\psi_{k}\left(  \xi^{\left(
q\right)  }\right)  \,w^{\left(  q\right)  },
\]
where~$\xi^{\left(  q\right)  }$ are collocation (quadrature) points,
and~$w^{\left(  q\right)  }$ are quadrature weights. We refer, e.g.,
to~\cite{LeMaitre-2010-SMU} for an overview and discussion of integration rules.

\begin{figure}[ptbh]
\begin{center}
\includegraphics[width=12.8cm]{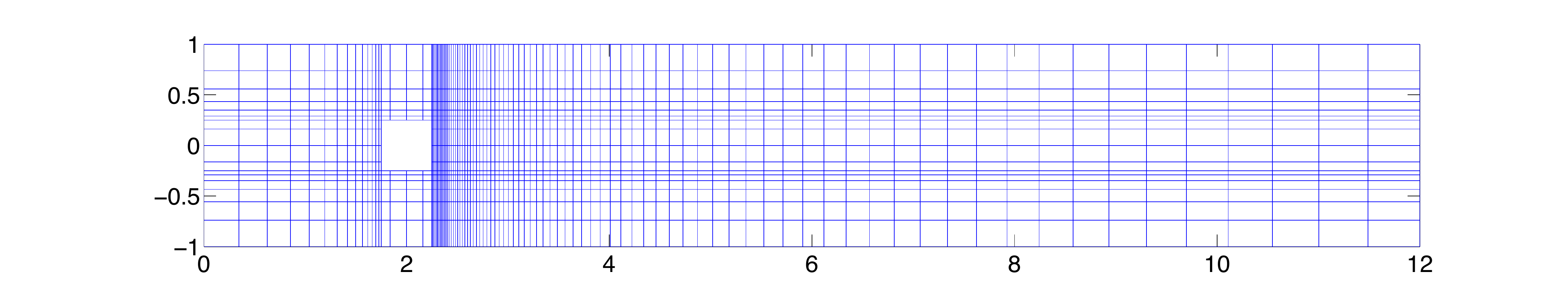}
\end{center}
\caption{Finite element mesh for the flow around an obstacle problem. }%
\label{fig:mesh-obstacle}%
\end{figure}

\subsection{Example: flow around an obstacle}

\label{sec:obstacle} In this section, we present results of numerical
experiments for a model problem given by a flow around an obstacle in a
channel of length $12$ and height $2$. 
The viscosity (\ref{eq:viscosity-gPC}) was taken to be a truncated lognormal
process with mean values $\nu_{0}=1/50$ or $1/150$, which corresponds to mean
Reynolds numbers $\operatorname{Re}_{0}=100$ or $300$, respectively, and its
representation was computed from an underlying Gaussian random process using
the transformation described in \cite{Ghanem-1999-NGS}. 
That is, for $\ell=0,\dots,M_{\nu}-1$, $\psi_{\ell}\left(  \xi\right)  $ is the product of
$N$ univariate Hermite polynomials, and denoting the coefficients of the
Karhunen-Lo\`{e}ve expansion of the Gaussian process by~$g_{j}\left(
x\right)  $ and $\eta_{j}=\xi_{j}-g_{j}$, $j=1,\dots,N$, the coefficients in
the expansion~(\ref{eq:viscosity-gPC}) are computed as
\[
\nu_{\ell}\left(  x\right)  =\mathbb{E}\left[  \psi_{\ell}\left(  \eta\right)
\right]  \exp\left[  g_{0}\left(  x\right)  +\frac{1}{2}\sum_{j=1}^{N}\left(
g_{j}\left(  x\right)  \right)  ^{2}\right]  .
\]
The covariance function of the Gaussian field, for points~$X_{1}=(x_{1}%
,y_{1})$, $X_{2}=(x_{2},y_{2})\in D$, was chosen to be
\[
C\left(  X_{1},X_{2}\right)  =\sigma_{g}^{2}\exp\left(  -\frac{\left\vert
x_{2}-x_{1}\right\vert }{L_{x}}-\frac{\left\vert y_{2}-y_{1}\right\vert
}{L_{y}}\right)  ,
\]
where~$L_{x}$ and $L_{y}$\ are the correlation lengths of the random
variables~$\xi_{i}$, $i=1,\dots,N$, in the $x$ and $y$ directions,
respectively, and $\sigma_{g}$ is the standard deviation of the Gaussian
random field. The correlation lengths were set to be equal to $25\%$ of the
width and height of the domain, i.e. $L_{x}=3$ and $L_{y}=0.5$. The
coefficient of variation of the lognormal field, defined as $CoV=\sigma_{\nu
}/\nu_{0}$ where $\sigma_{\nu}$ is the standard deviation, was $10\%$ or
$30\%$. The stochastic dimension was $N=2$.
\textcolor{black}{The lognormal formulation was chosen to enable exploration of a general
random field in which the viscosity guaranteed to be positive.
See \cite{Bab-Nob-Temp} for an example of the use  of this formulation for diffusion problems
and \cite{Zhang:book} for its use in models of porous media.}

We implemented the methods in
\texttt{Matlab} using \texttt{IFISS~3.3}~\cite{ERS-SIREV}. 
The spatial discretization uses a stretched grid, discretized by $1520$ Taylor-Hood
finite\ elements; the domain and grid are shown in
Figure~\ref{fig:mesh-obstacle}.
There are $12,640$ velocity and $1640$ pressure degrees of freedom. 
the degree used for the polynomial
expansion of the solution was $P=3$, and the degree used for the expansion of
the lognormal process was $2P$, which ensures a complete representation of the
process in the discrete problem~\cite{Matthies-2005-GML}. With these settings,
$M=10$ and $M_{\nu}=\widehat{M}=28$, and $\mathbf{H}_{\ell}$ is of order $10$
in~(\ref{eq:linearized-system}).

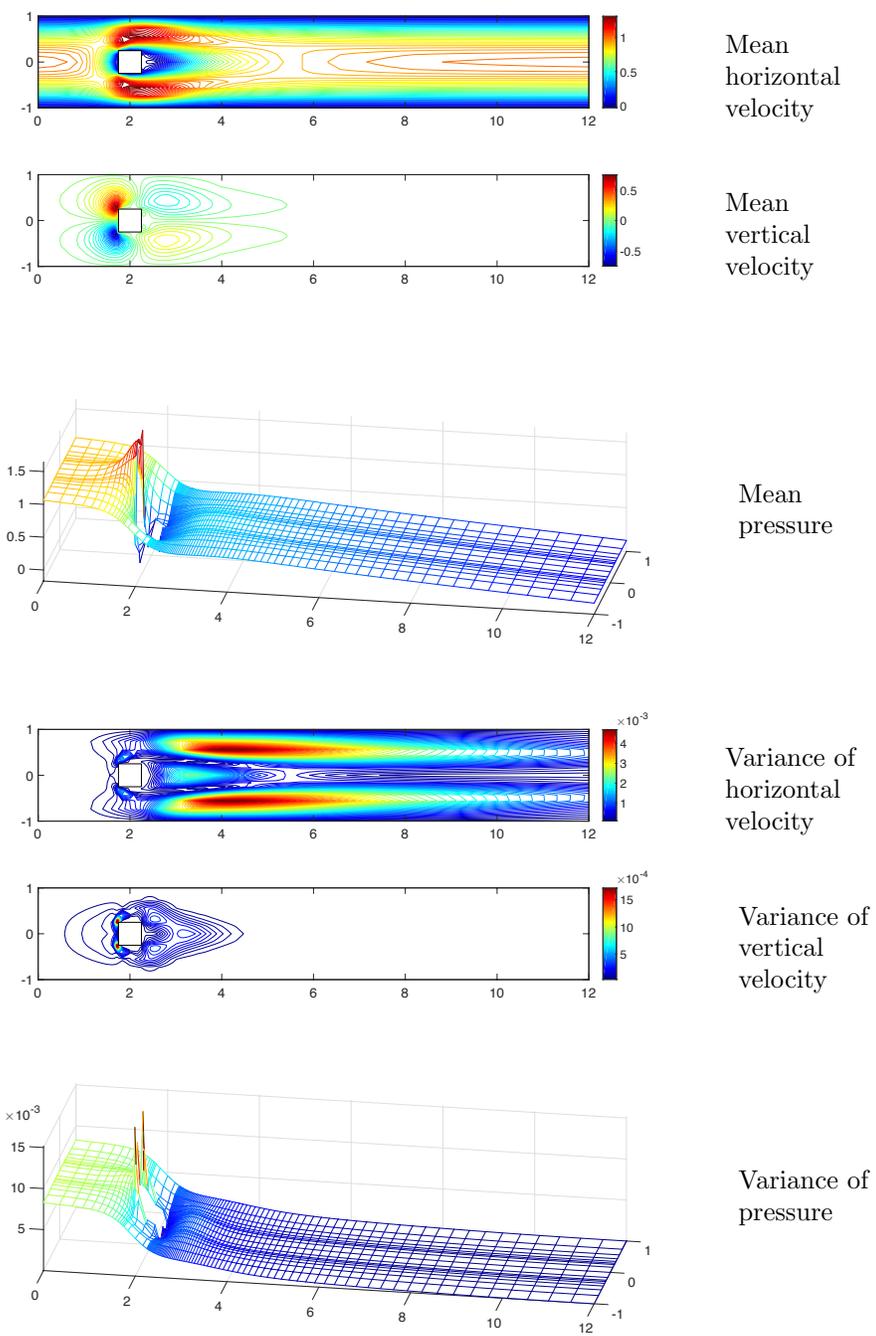
\begin{figure}[ptbh]
\begin{center}
\begin{picture}(500,515) (20,20)
\put(30,485) {\includegraphics[width=10cm]{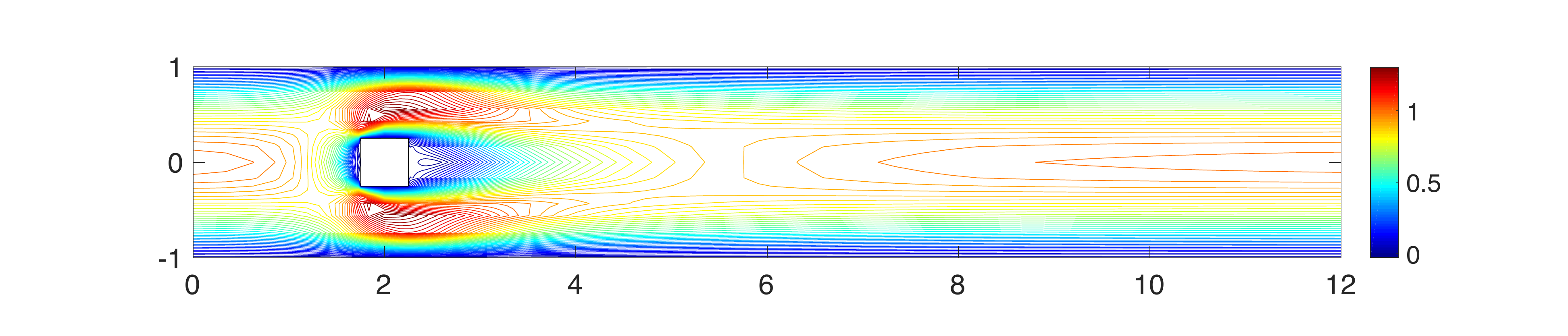}}
\put(325,520) {Mean} \put(325,508){horizontal}\put(325,496){velocity}
\put(30,425) {\includegraphics[width=10cm]{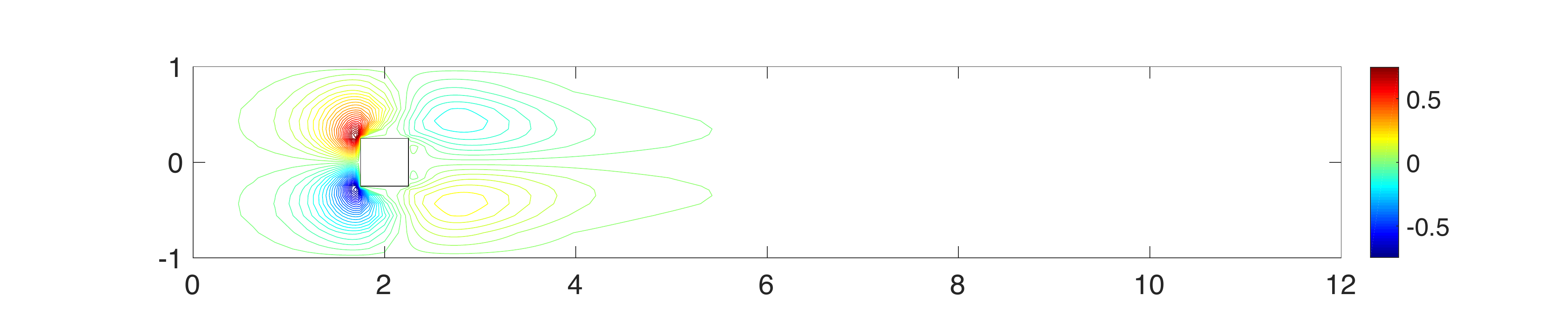}}
\put(325,460) {Mean} \put(325,448){vertical}\put(325,436){velocity}
\put(30,285) {\includegraphics[width=10cm]{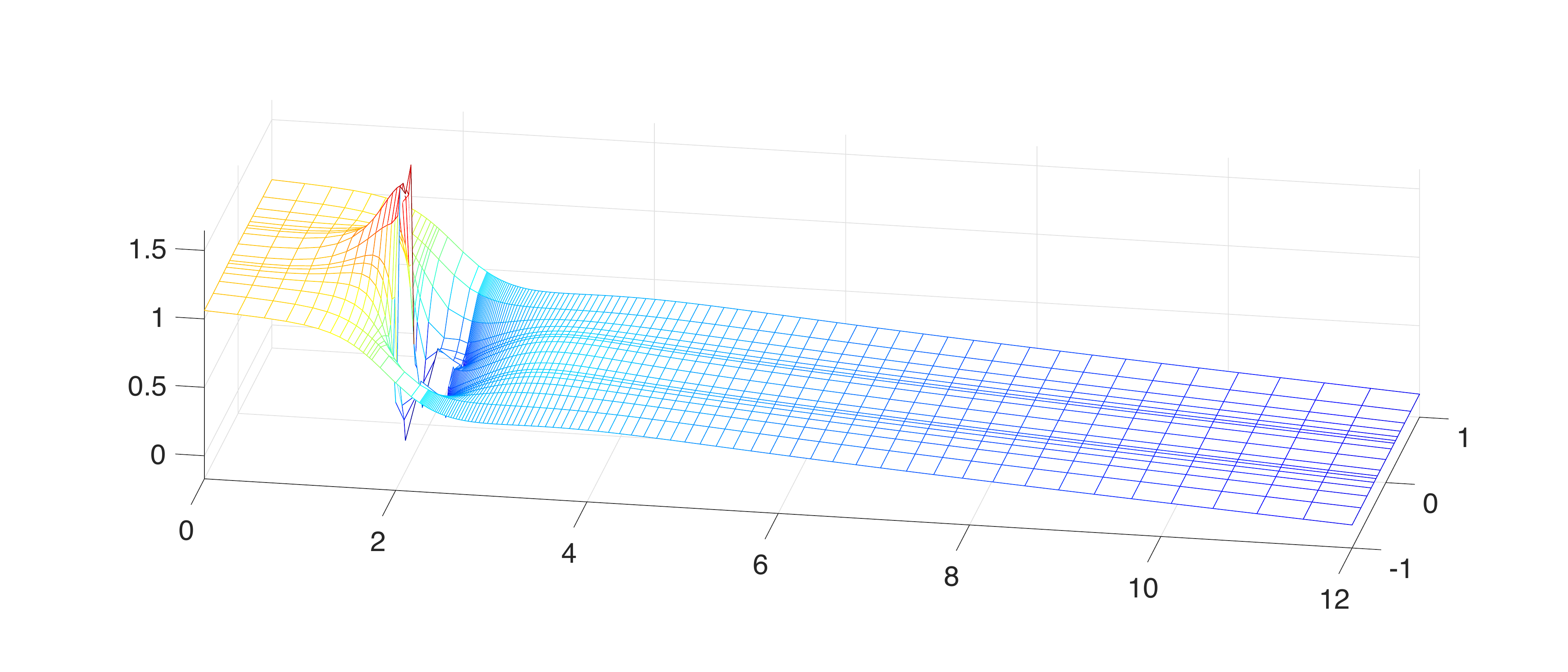}}
\put(330,350) {Mean} \put(330,338){pressure}
\put(30,215) {\includegraphics[width=10.0cm]{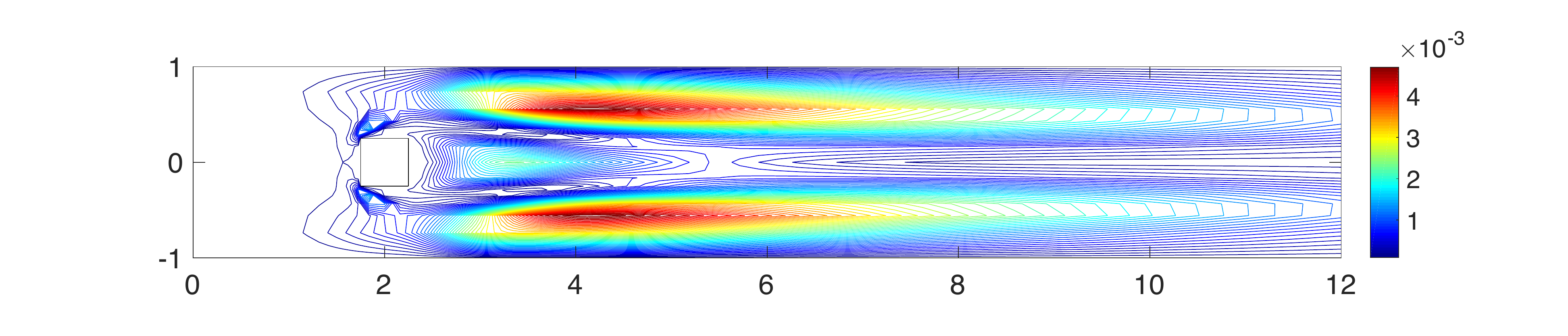}}
\put(325,250) {Variance of} \put(325,238){horizontal}\put(325,226){velocity}
\put(30,155) {\includegraphics[width=10.0cm]{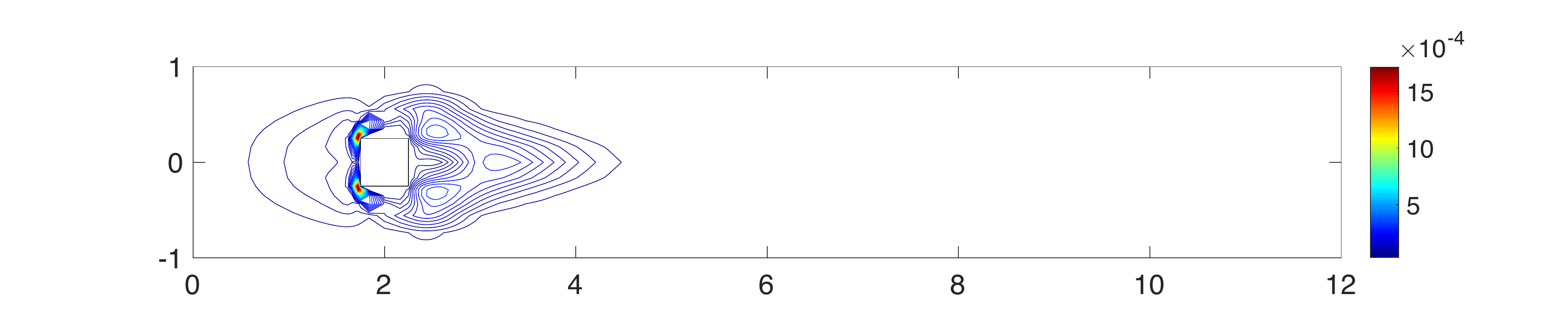}}
\put(330,190) {Variance of} \put(330,178){vertical} \put(330,166){velocity}
\put(30, 25) {\includegraphics[width=10.0cm]{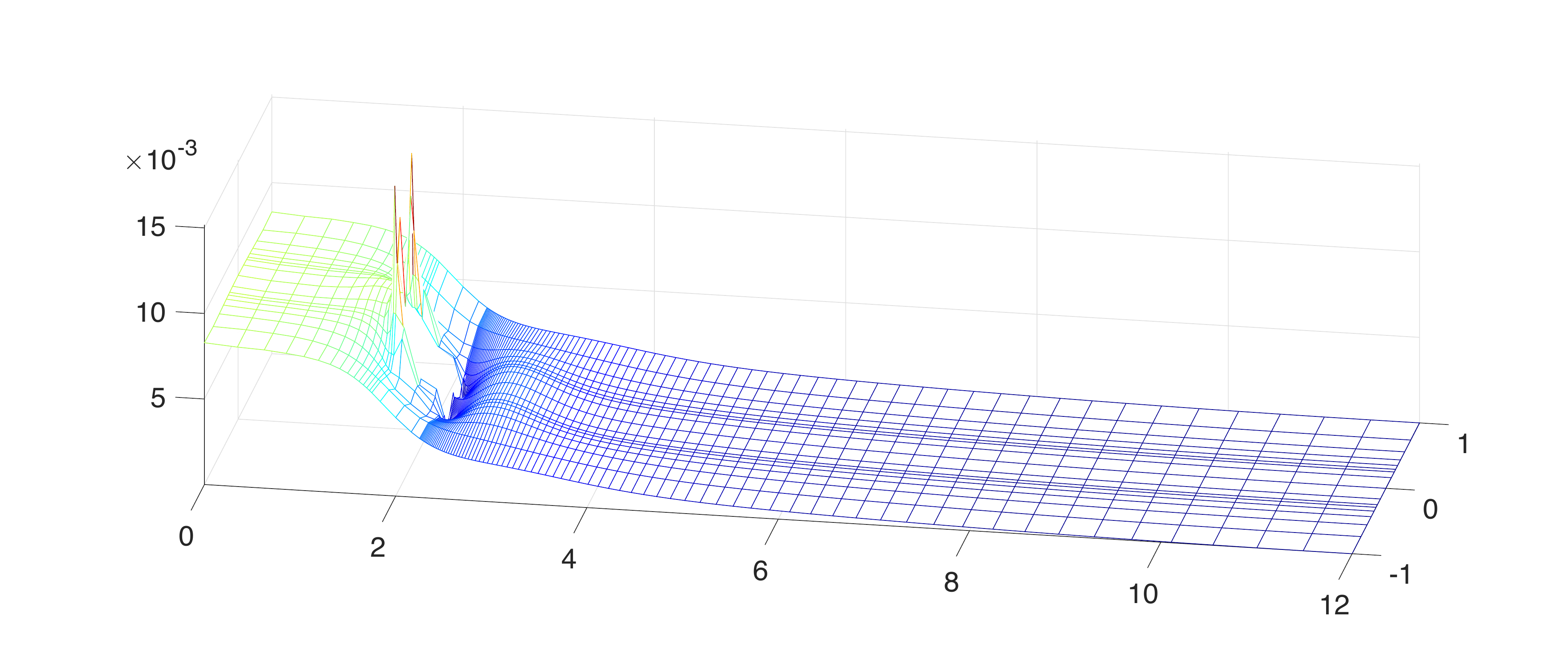}}
\put(330,90) {Variance of} \put(330,78){pressure}
\end{picture}
\end{center}
\caption{Mean horizontal and vertical velocities and pressure (top) and
variances of the same quantities (bottom), for $\operatorname{Re}=100$ and
$CoV=10\%$. }%
\label{fig:sRe100_CoV10_0}%
\end{figure}

\begin{figure}[ptbh]
\begin{center}
\includegraphics[width=10.0cm]{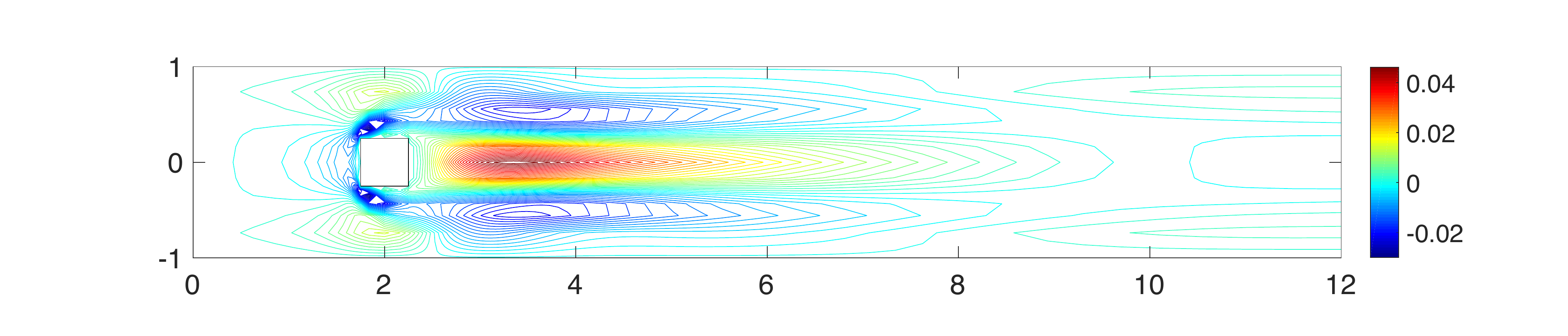}
\includegraphics[width=10.0cm]{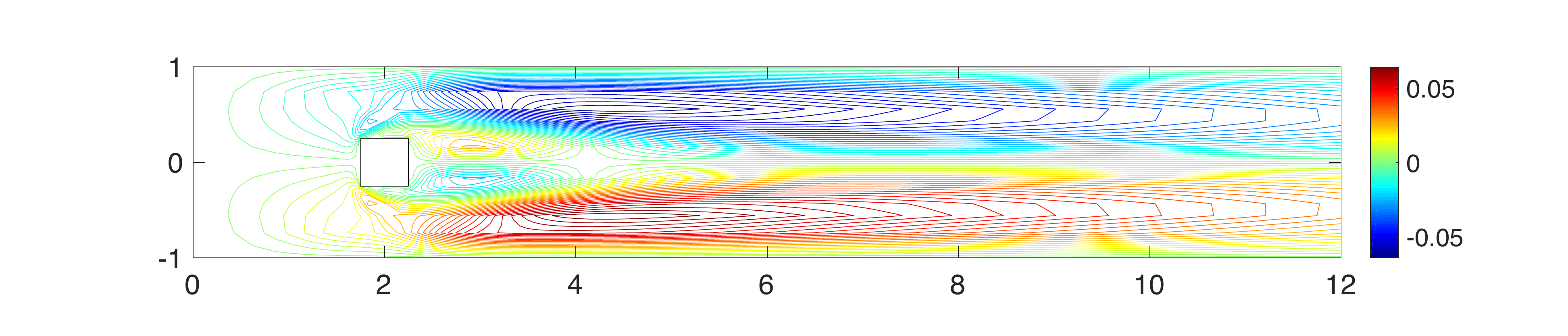}
\includegraphics[width=10.0cm]{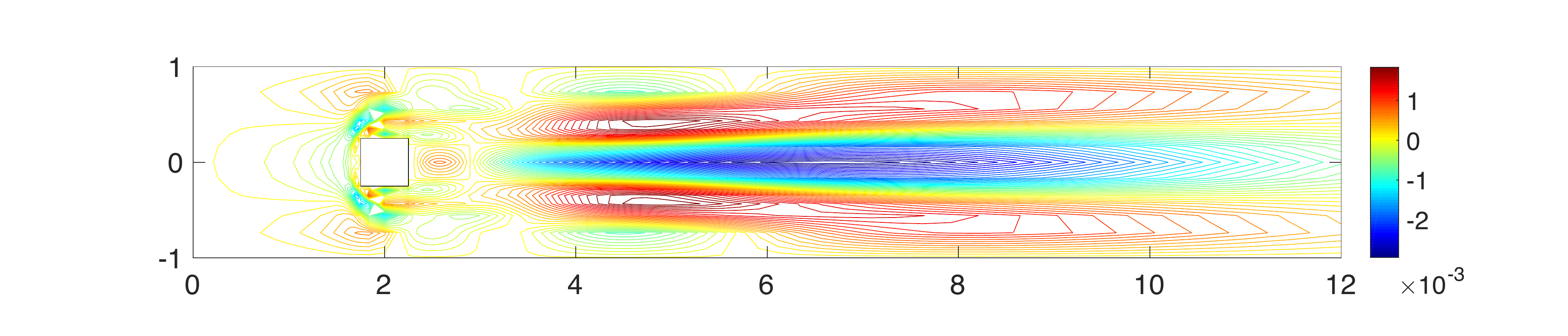}
\includegraphics[width=10.0cm]{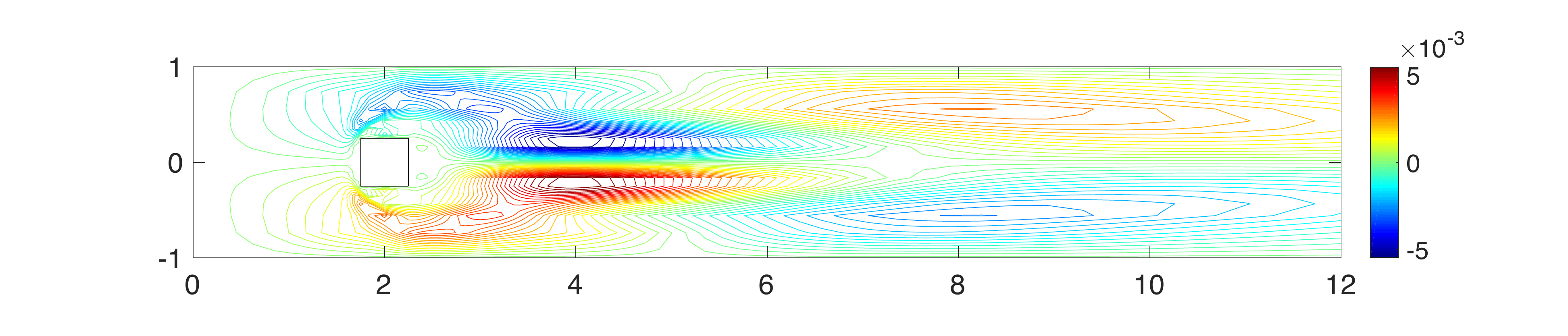}
\end{center}
\caption{Coefficients~$1-4$ of the gPC expansion of the horizontal velocity,
$\operatorname{Re}_{0}=100$ and $CoV=10\%$.}%
\label{fig:sRe100_CoV10_gPC}%
\end{figure}

Consider first the case of $\operatorname{Re}_{0}=100$ and $CoV=10\%$.
Figure~\ref{fig:sRe100_CoV10_0} shows the mean horizontal and vertical
components of the velocity and the mean pressure (top), and the variances of
the same quantities (bottom). It can be seen that there is symmetry in all the
quantities, the mean values are essentially the same as we would expect in the
deterministic case, and the variance of the horizontal velocity component is
concentrated in two \textquotedblleft eddies\textquotedblright\ and is larger
than the variance of the vertical velocity component.
Figure~\ref{fig:sRe100_CoV10_gPC} illustrates values of several coefficients
of expansion~(\ref{eq:u-gPC}) of the horizontal velocity. All the coefficients
are symmetric, and as the index increases they become more oscillatory and
their values decay. We found the same trends for the coefficients of the
vertical velocity component and of the pressure. Our observations are
qualitatively consistent with numerical experiments of Powell and
Silvester~\cite{Powell-2012-PSS}.

\begin{figure}[ptbh]
\begin{center}
\begin{picture}(500,235)(20,0)
\put(30,185) {\includegraphics[width=10.0cm]{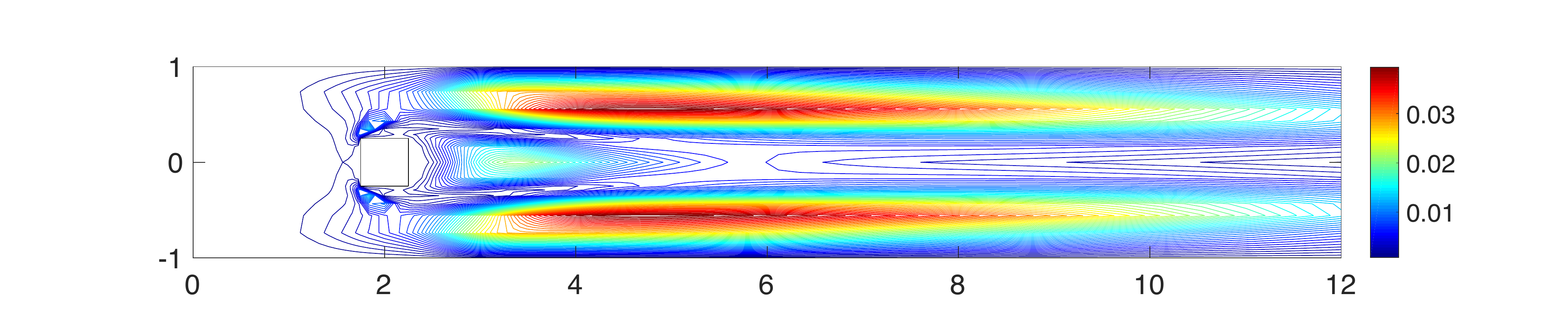}}
\put(325,220) {Variance of} \put(325,208){horizontal}\put(325,196){velocity}
\put(30,125) {\includegraphics[width=10.0cm]{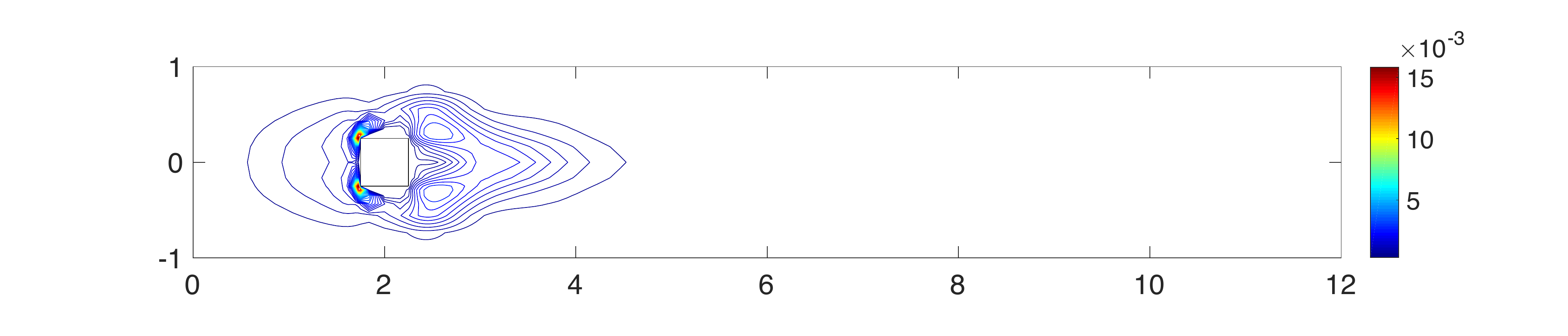}}
\put(330,160) {Variance of} \put(330,148){vertical} \put(330,136){velocity}
\put(30, -5) {\includegraphics[width=10.0cm]{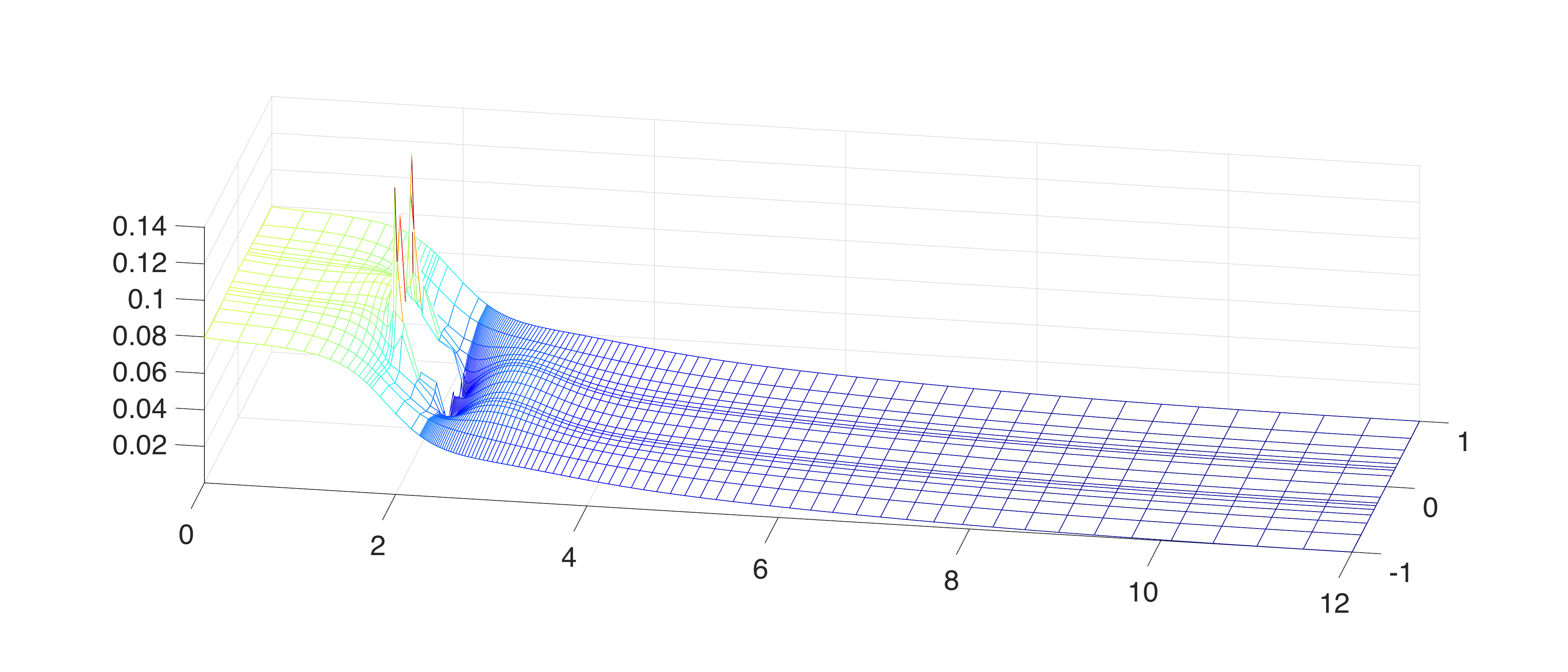}}
\put(330,50) {Variance of} \put(330,38){pressure}
\end{picture}
\end{center}
\caption{Variances of velocity components and pressure for $\operatorname{Re}%
_{0}=100$ and $CoV=30\%$.}%
\label{fig:sRe100_CoV30_var}%
\end{figure}
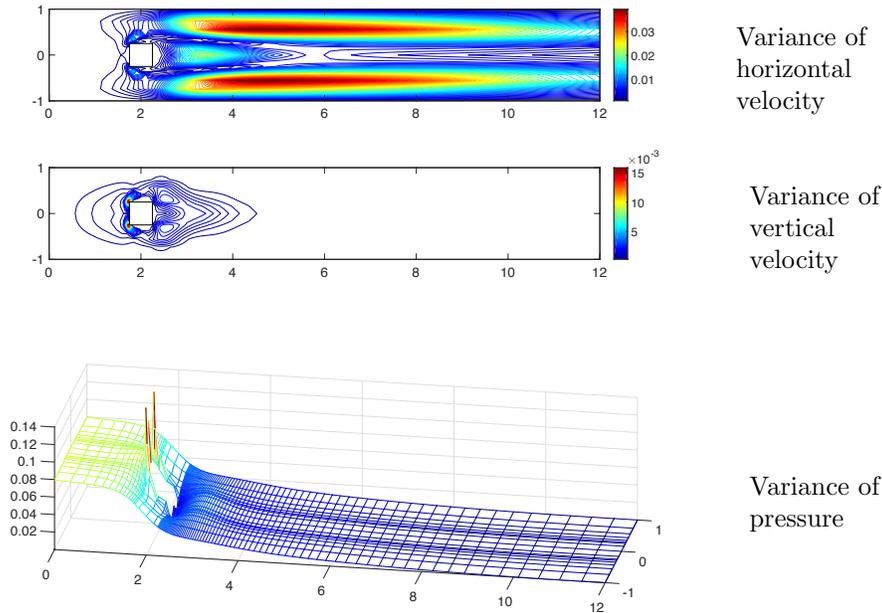

We also tested (the same) $\operatorname{Re}_{0}=100$ with increased
$CoV=30\%$. We found that the mean values are essentially the same as in the
previous case; Figure~\ref{fig:sRe100_CoV30_var} shows the variances, which
display the same qualitative behavior but have values that are approximately
$10$ times larger than for the case $CoV=10\%$.

\begin{figure}[t]
\begin{center}
\includegraphics[width=6.0cm]{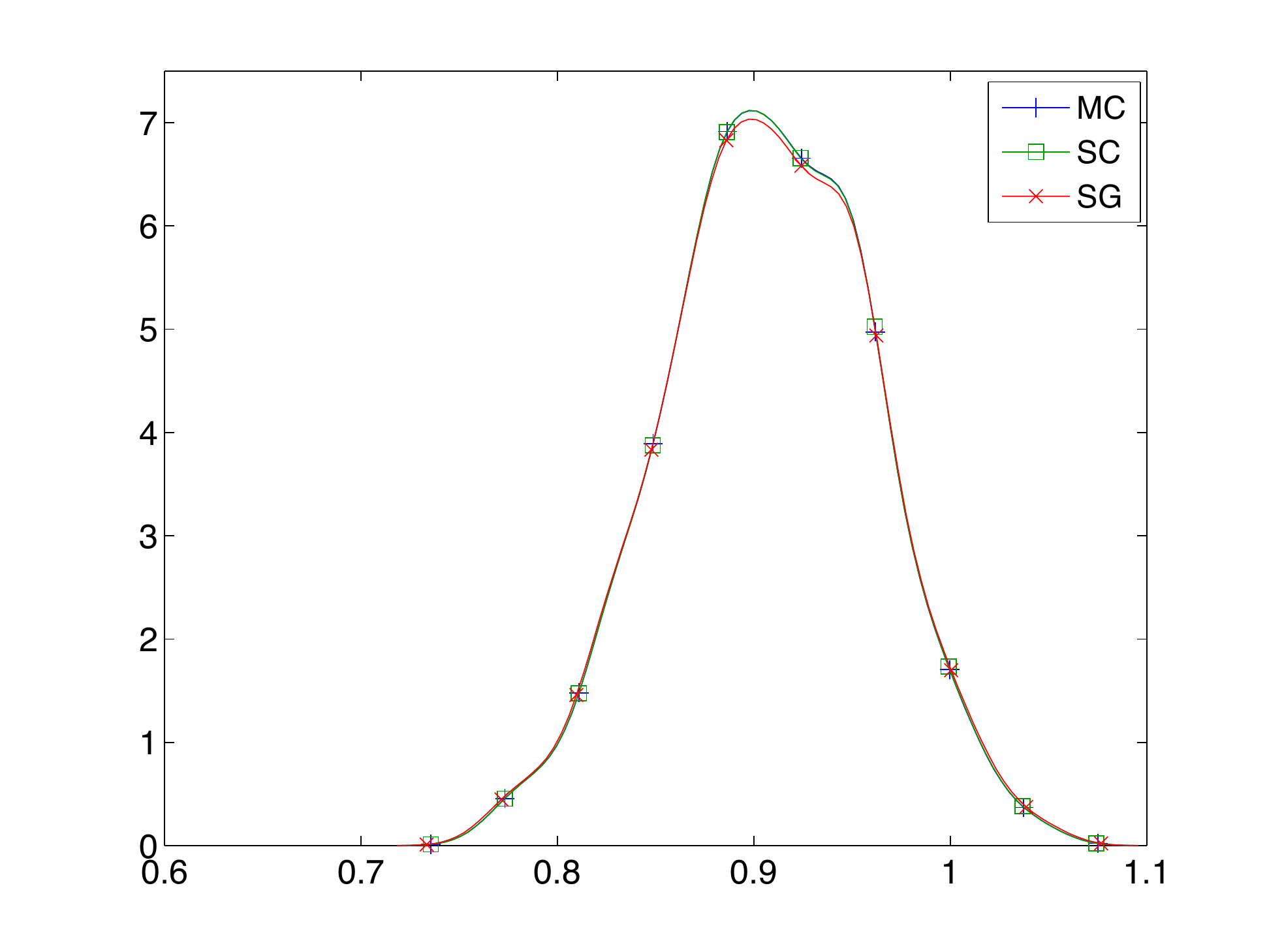}
\includegraphics[width=6.0cm]{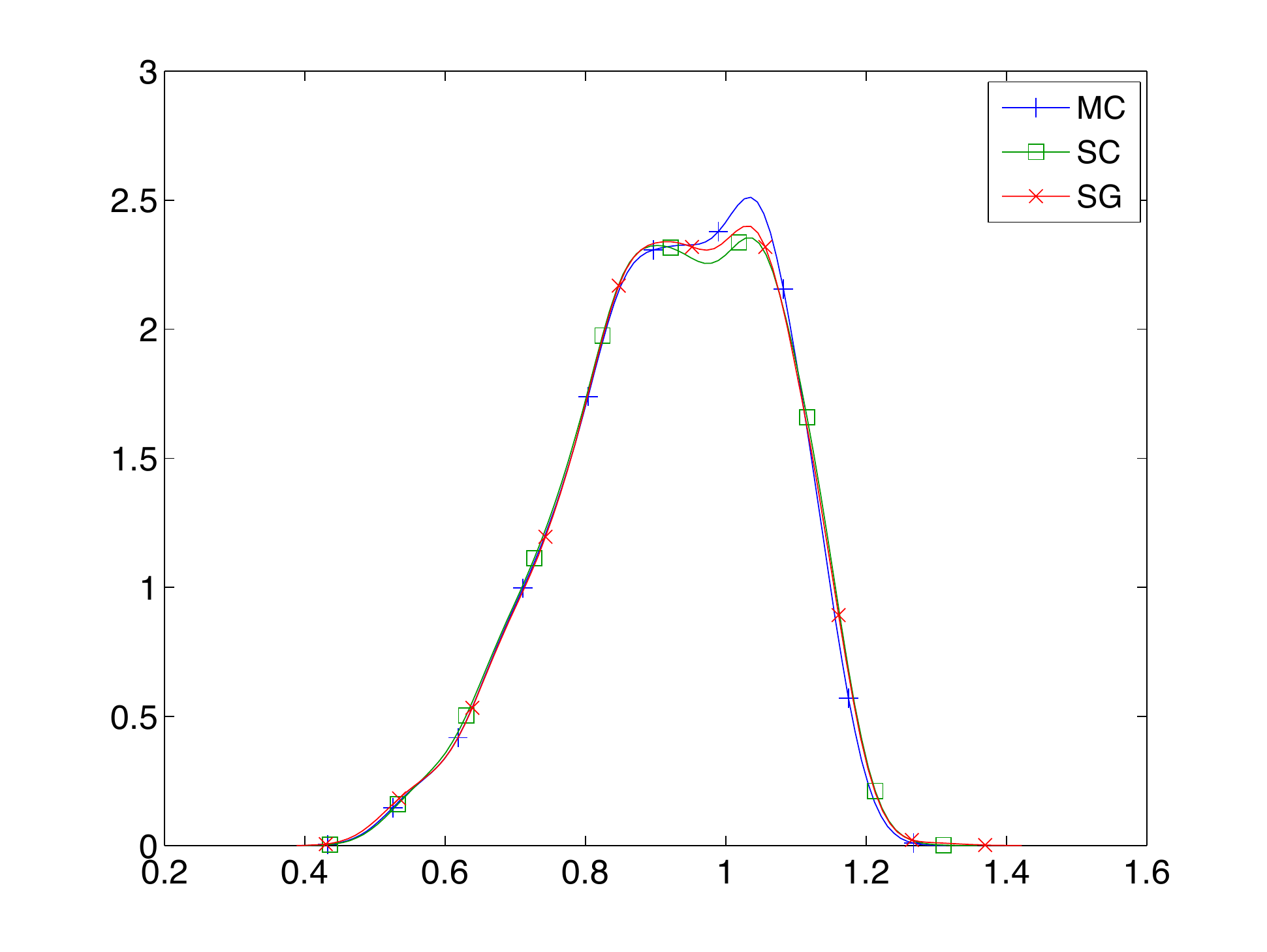}
\newline\includegraphics[width=6.0cm]{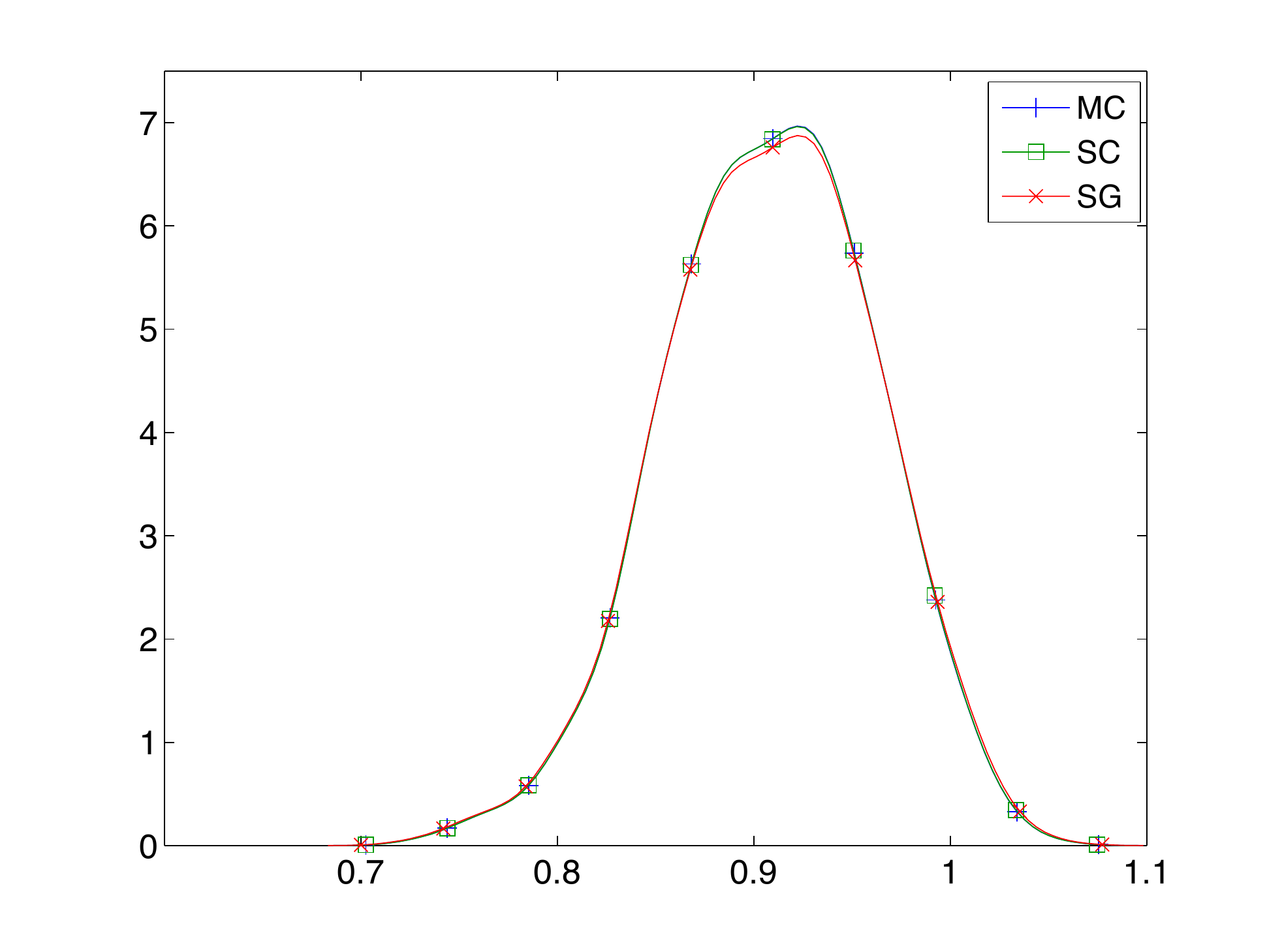}
\includegraphics[width=6.0cm]{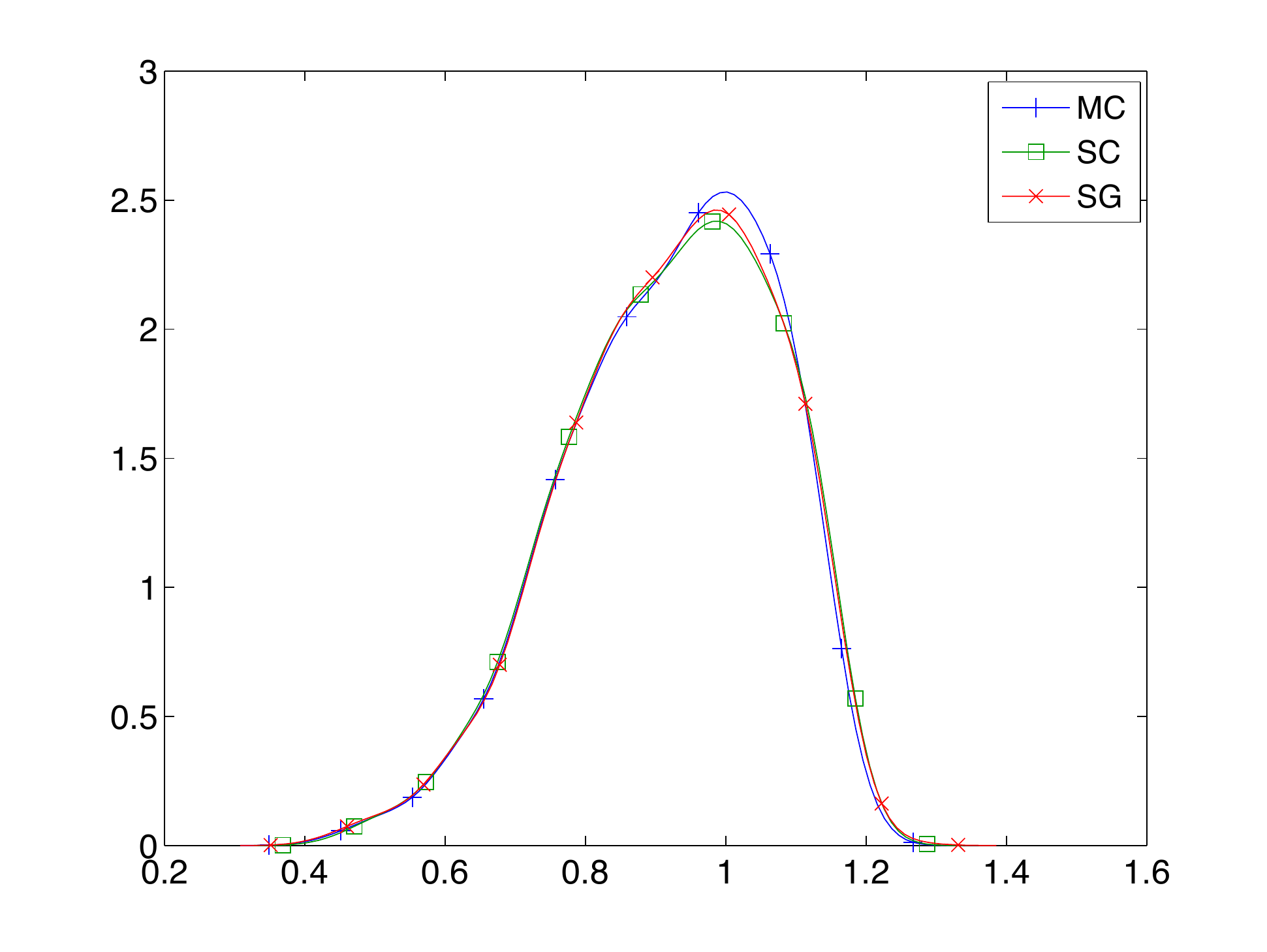}
\end{center}
\caption{Estimated pdfs of the velocities $u_{x}$ with $\operatorname{Re}%
_{0}=100$, $CoV=10\%$ (left) and $30\%$ (right) at the points with coordinates
$(4.0100,-0.4339)$ (top) and $(4.0100,0.4339)$ (bottom). }%
\label{fig:sRe100_QoI}%
\end{figure}

A different perspective on the results is given in Figure~\ref{fig:sRe100_QoI}%
, which shows estimates of the probability density function (pdf) for the
horizontal velocity at two points in the domain, $(4.0100,-0.4339)$ and
$(4.0100,0.4339)$. These are locations at which large variances of the
solution were seen, see Figures~\ref{fig:sRe100_CoV10_0}
and~\ref{fig:sRe100_CoV30_var}. The results were obtained using
\texttt{Matlab}'s \texttt{ksdensity} function. It can be seen that with the
larger value of $CoV$, the support of the velocity pdf is wider, and except
for the peak values, for fixed $CoV$ the shapes of the pdfs at the two points
are similar, indicating a possible symmetry of the stochastic solution. For
this benchmark, we also obtained analogous data using the Monte Carlo and
collocation sampling methods; it can be seen from the figure that these
methods produced similar results. \footnote{The results for Monte Carlo were
obtained using $10^{3}$ samples, and those for collocation were found using a
Smolyak sparse grid with Gauss-Hermite quadrature and grid level~$n_{q}=4$.}

\begin{figure}[ptbh]
\begin{center}
\begin{picture}(500,515) (20,20)
\put(30,485) {\includegraphics[width=10cm]{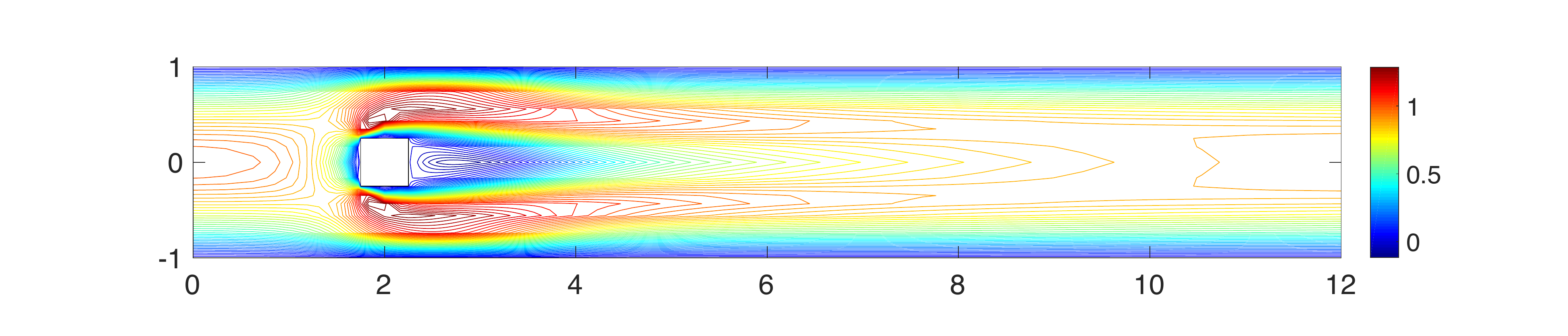}}
\put(325,520) {Mean} \put(325,508){horizontal}\put(325,496){velocity}
\put(30,425) {\includegraphics[width=10cm]{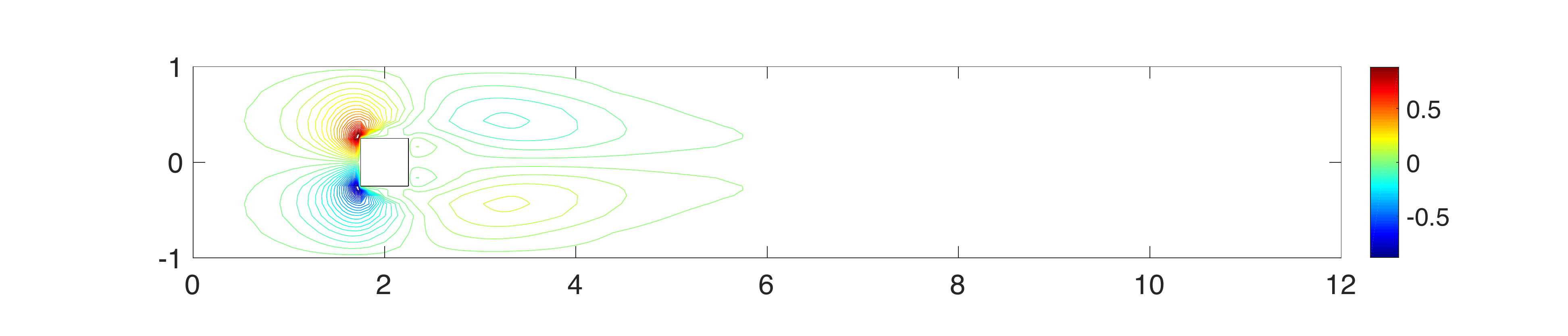}}
\put(325,460) {Mean} \put(325,448){vertical}\put(325,436){velocity}
\put(30,285) {\includegraphics[width=10cm]{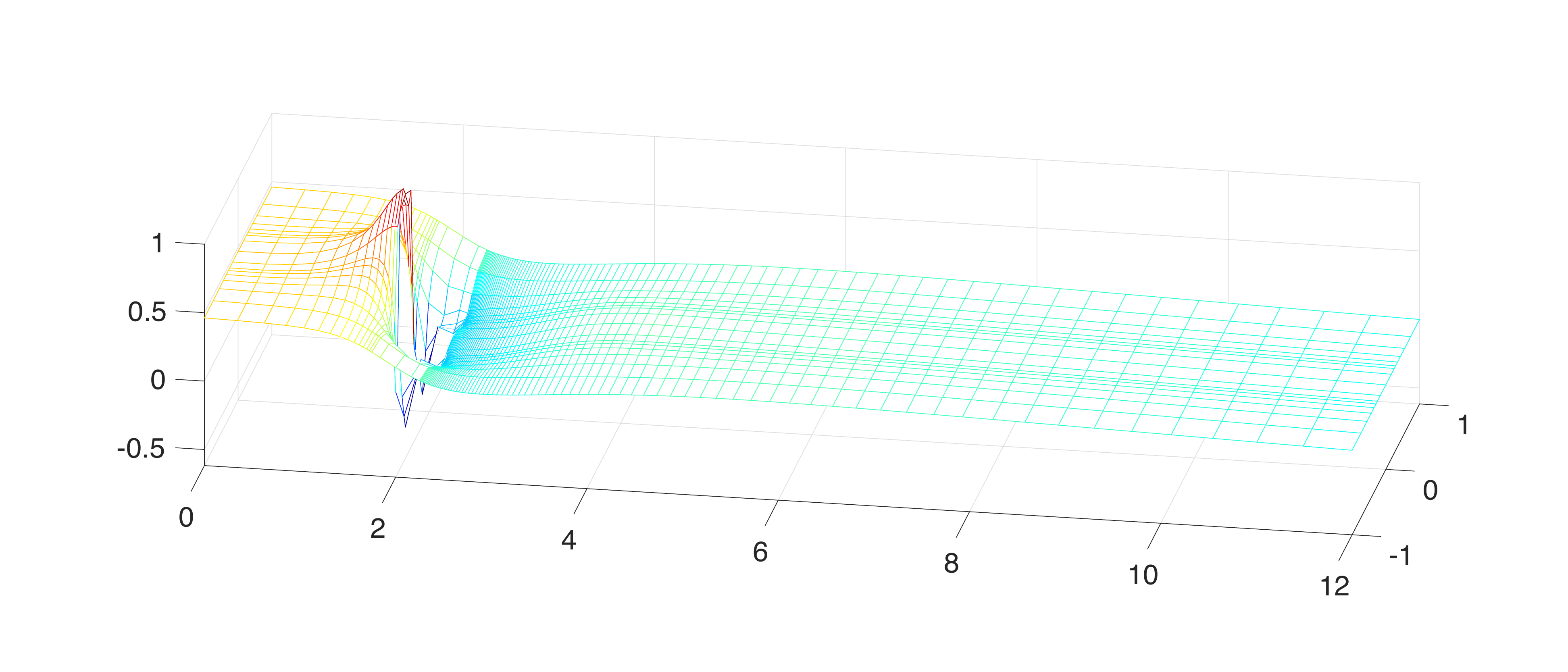}}
\put(330,350) {Mean} \put(330,338){pressure}
\put(30,215) {\includegraphics[width=10.0cm]{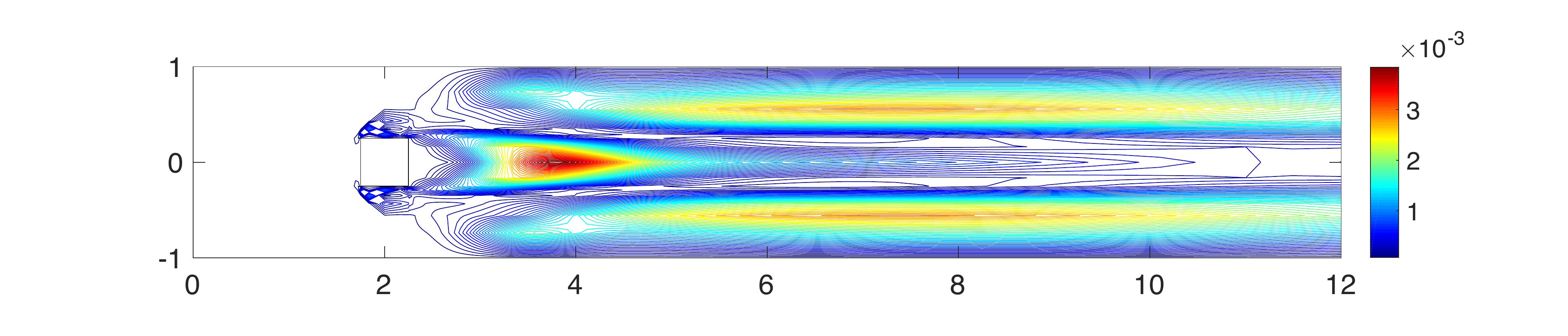}}
\put(325,250) {Variance of} \put(325,238){horizontal}\put(325,226){velocity}
\put(30,155) {\includegraphics[width=10.0cm]{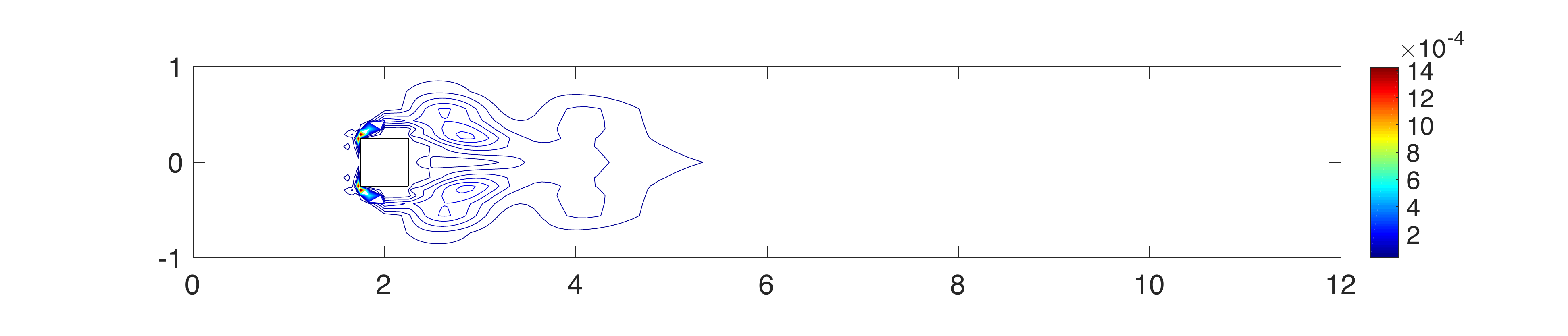}}
\put(330,190) {Variance of} \put(330,178){vertical} \put(330,166){velocity}
\put(30, 25) {\includegraphics[width=10.0cm]{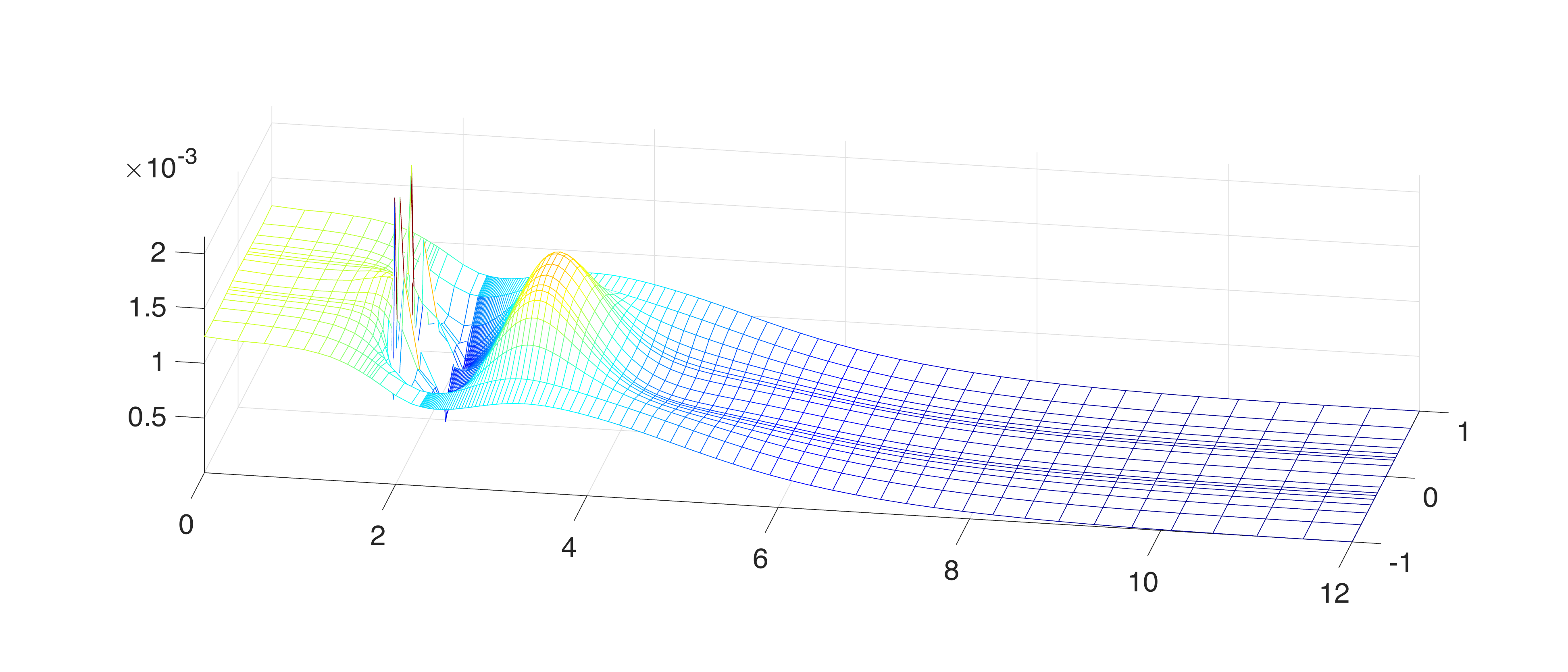}}
\put(330,90) {Variance of} \put(330,78){pressure}
\end{picture}
\end{center}
\caption{Mean horizontal and vertical velocities and pressure (top) and
variances of the same quantities (bottom), for $\operatorname{Re}_{0}=300$ and
$CoV=10\%$. }%
\label{fig:sRe300_CoV10_0}%
\end{figure}

Next, we consider a larger value of the mean Reynolds number,
$\operatorname{Re}_{0}=300$. Figure~\ref{fig:sRe300_CoV10_0} shows the means
and variances for the velocities and pressure for $CoV=10\%$. It is evident
that increased $\operatorname{Re}_{0}$ results in increased values of the mean
quantities, but they are again similar to what would be expected in the
deterministic case. The variances exhibit wider eddies than for
$\operatorname{Re}_{0}=100$, and in this case there is only one region of the
largest variance in the horizontal velocity, located just to the right of the
obstacle; this is also a region with increased variance of the pressure.

\begin{figure}[ptbh]
\begin{center}
\begin{picture}(500,235)(20,0)
\put(30,185) {\includegraphics[width=10.0cm]{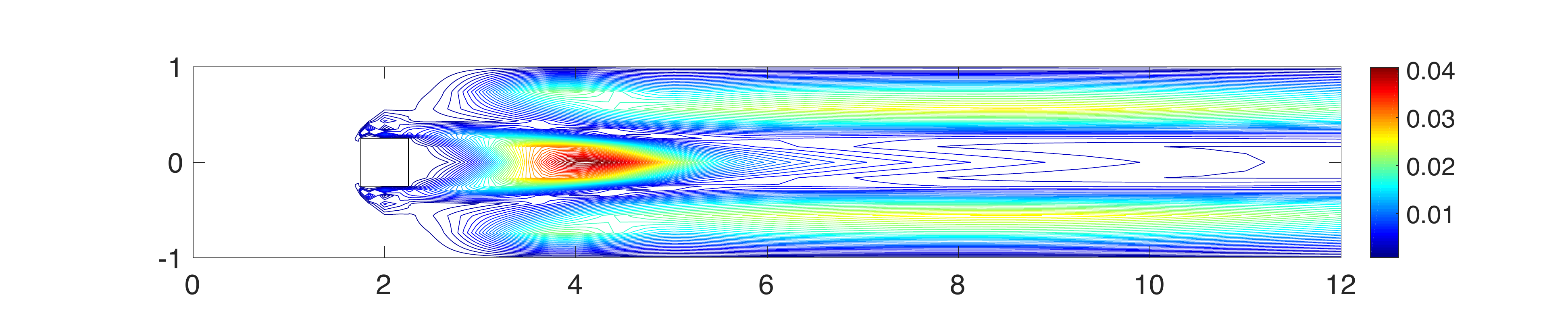}}
\put(325,220) {Variance of} \put(325,208){horizontal}\put(325,196){velocity}
\put(30,125) {\includegraphics[width=10.0cm]{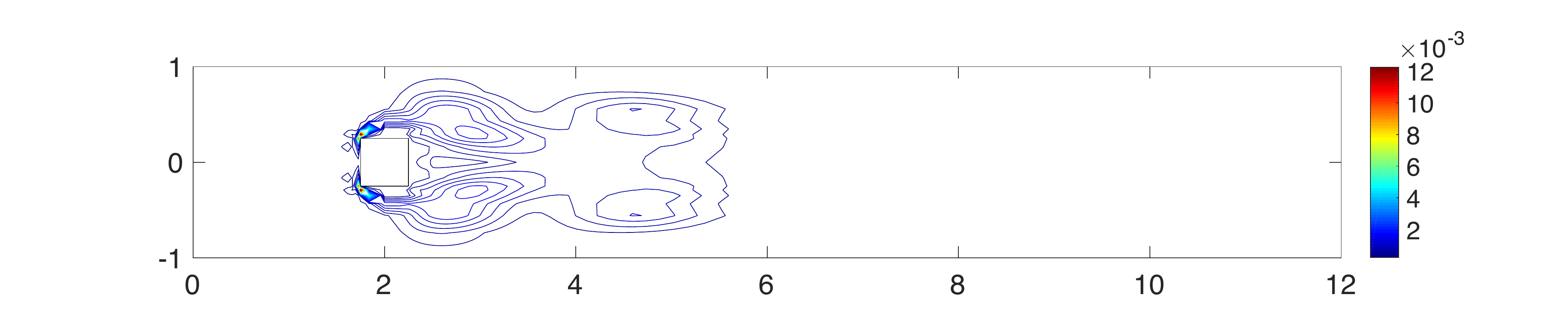}}
\put(330,160) {Variance of} \put(330,148){vertical} \put(330,136){velocity}
\put(30, -5) {\includegraphics[width=10.0cm]{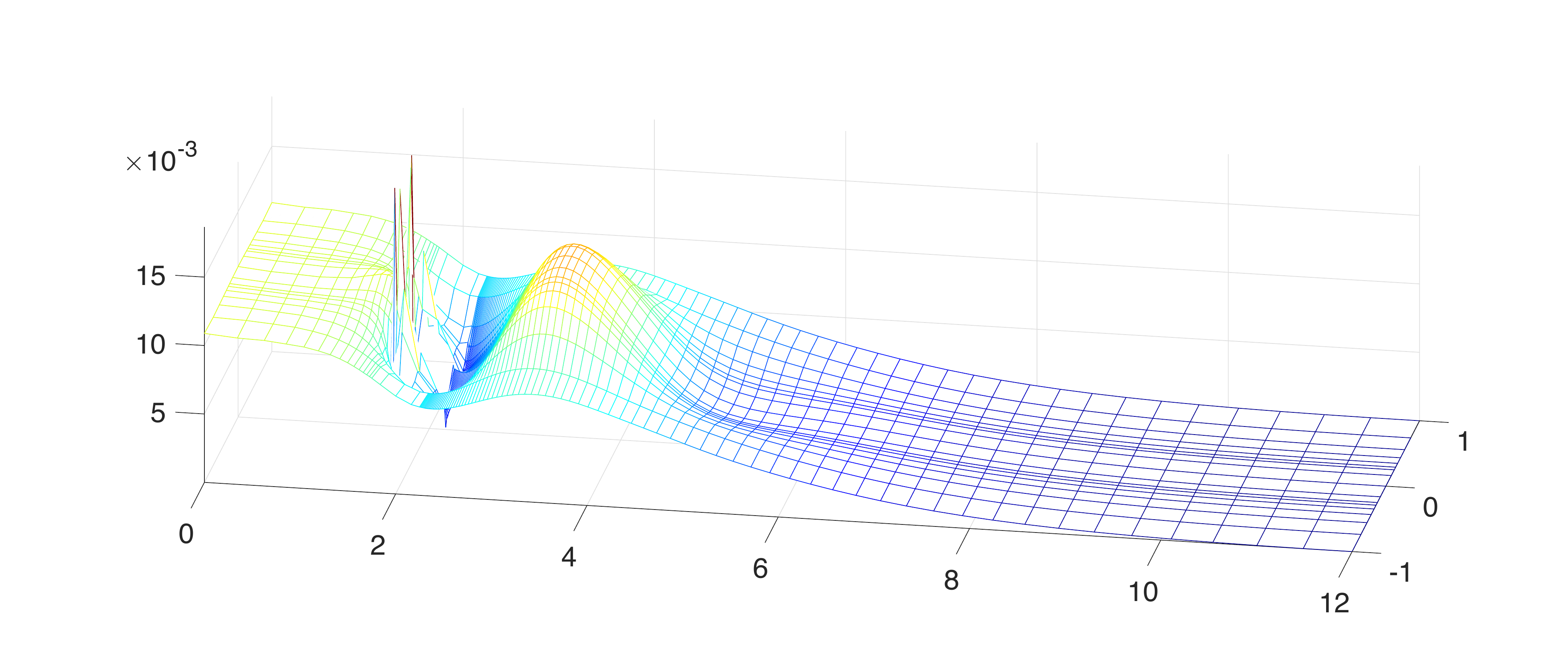}}
\put(330,50) {Variance of} \put(330,38){pressure}
\end{picture}
\end{center}
\caption{Variances of velocity components and pressure for $\operatorname{Re}%
_{0}=300$ and $CoV=30\%$.}%
\label{fig:sRe300_CoV30_var}%
\end{figure}

In similar tests for the larger value $CoV=30\%$, we found that the mean
values are essentially the same as for $CoV=10\%$, and
Figure~\ref{fig:sRe300_CoV30_var} shows the variances of velocities and
pressures. From the figure it can be seen that the variances are qualitatively
the same but approximately $10$ times larger than for $CoV=10\%$, results
similar to those found for $\operatorname{Re}_{0}=100$.

\begin{figure}[t]
\begin{center}
\includegraphics[width=6.0cm]{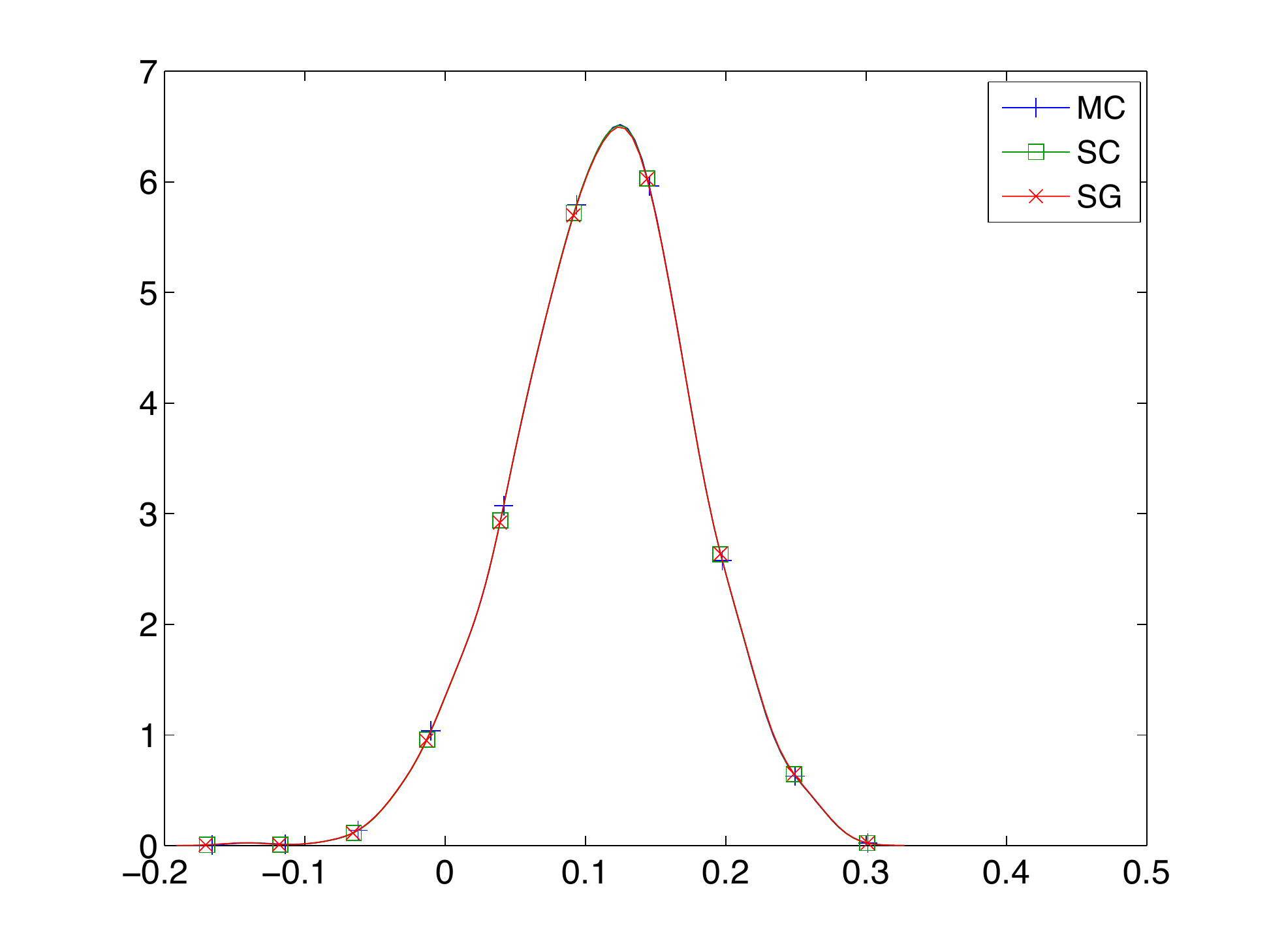}
\includegraphics[width=6.0cm]{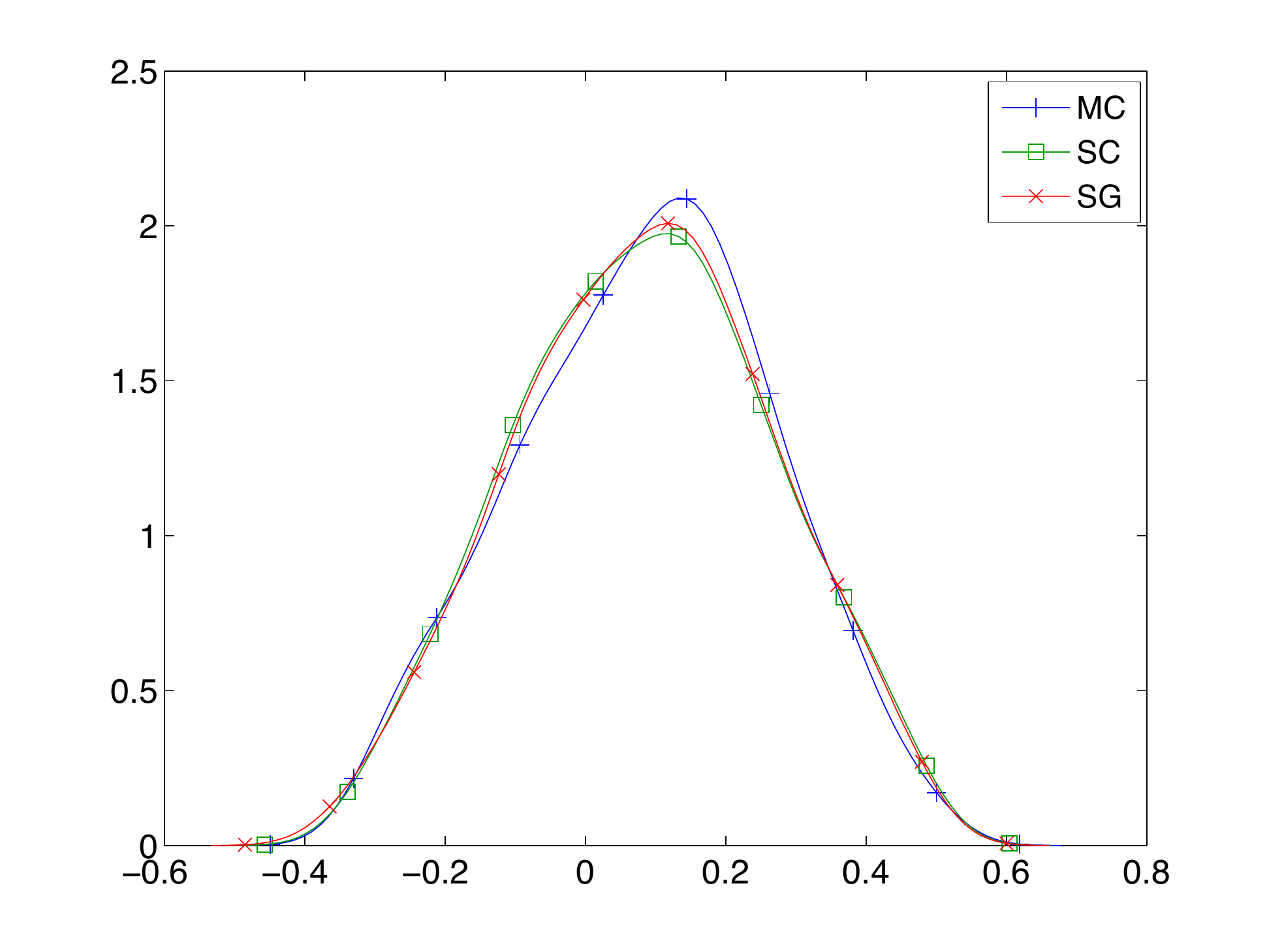}
\newline\includegraphics[width=6.0cm]{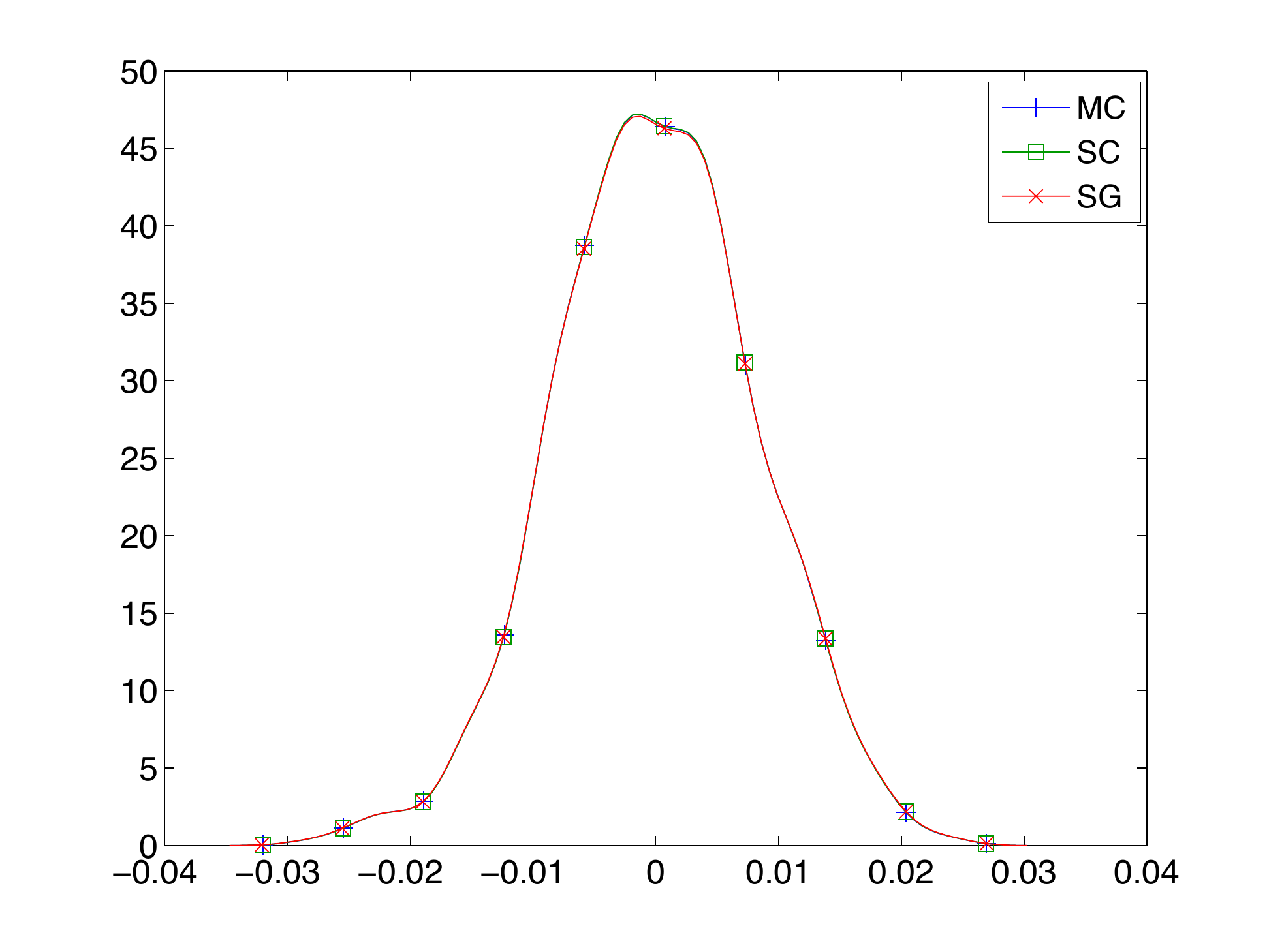}
\includegraphics[width=6.0cm]{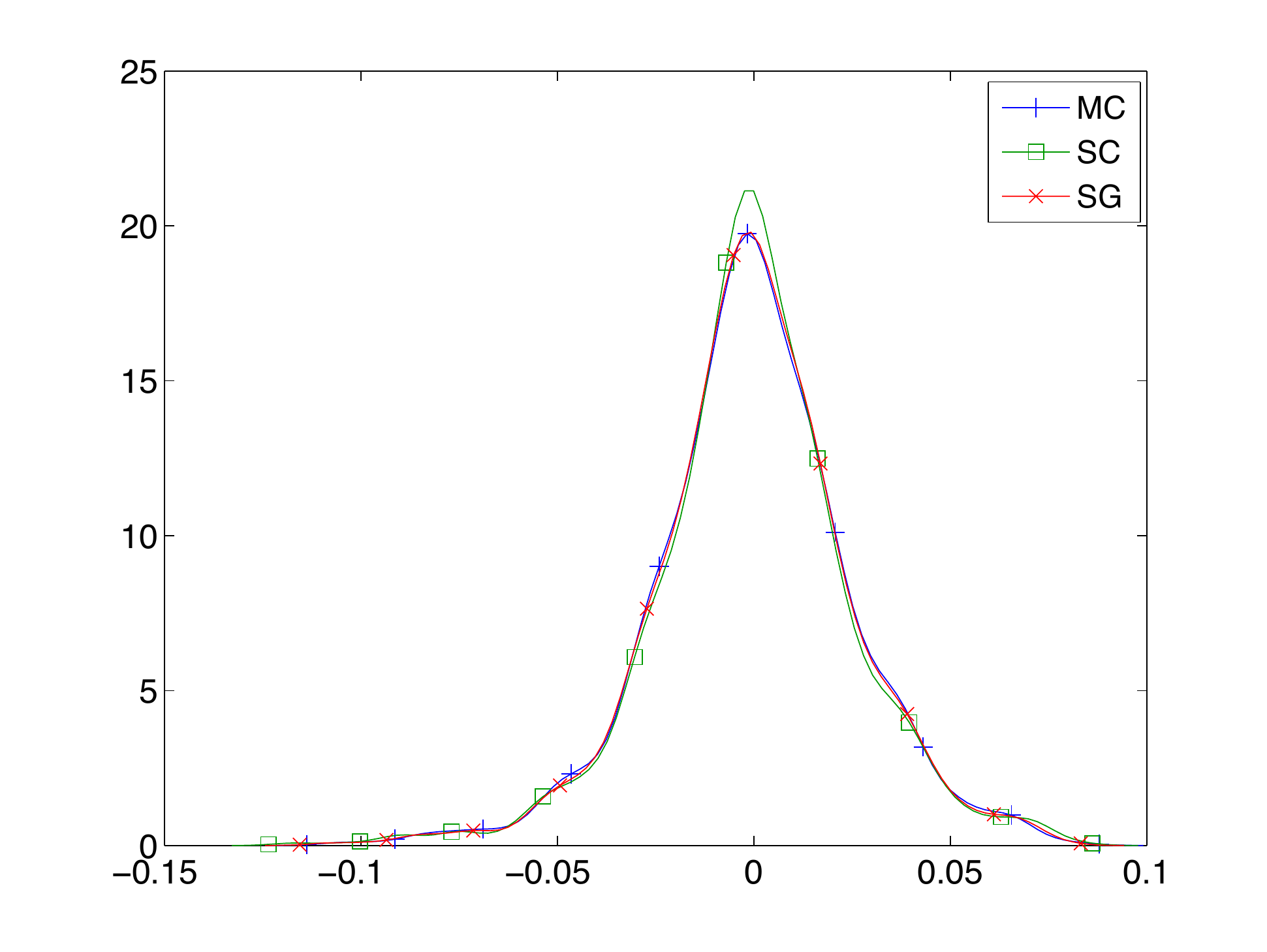}
\newline\includegraphics[width=6.0cm]{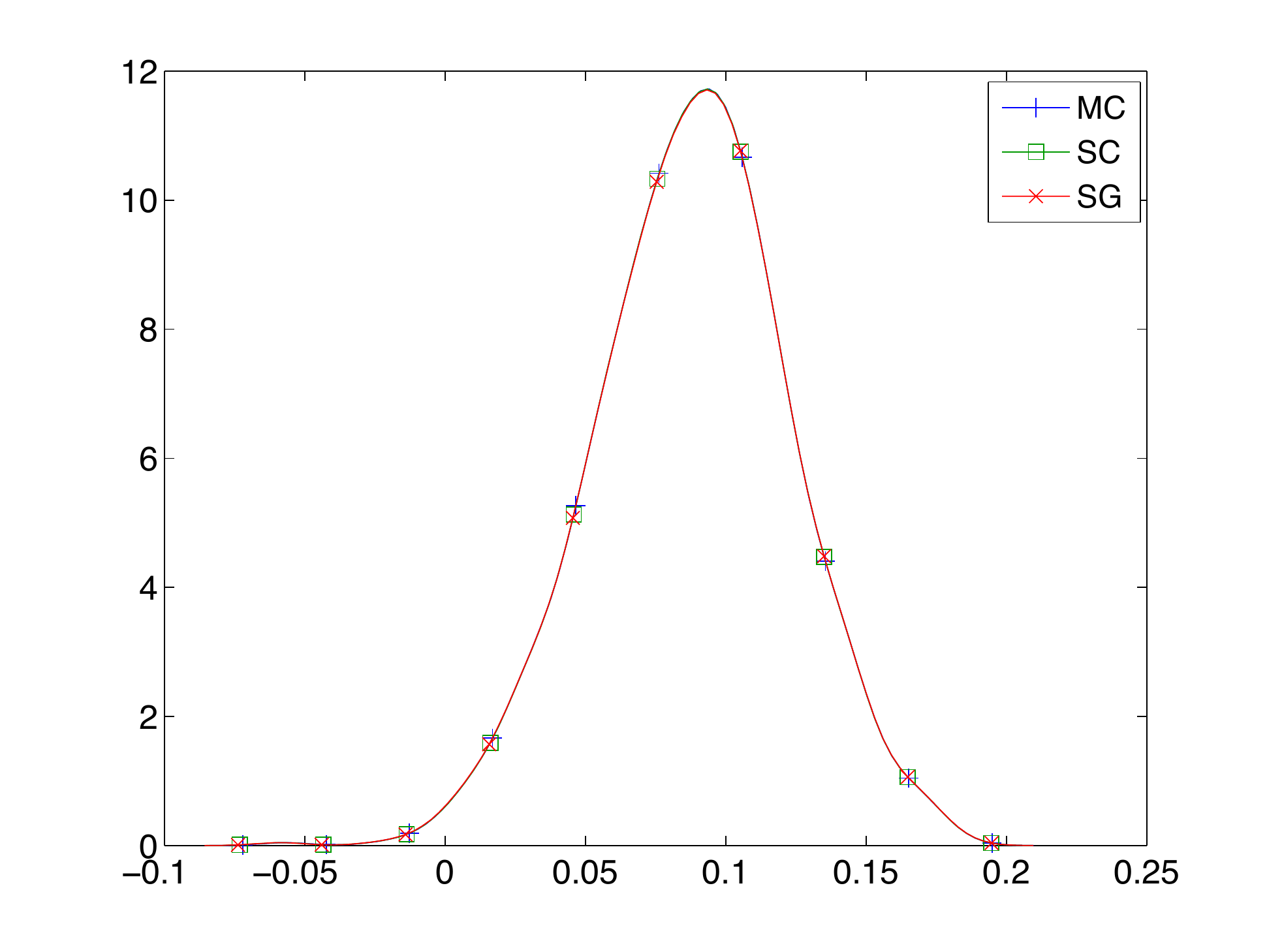}
\includegraphics[width=6.0cm]{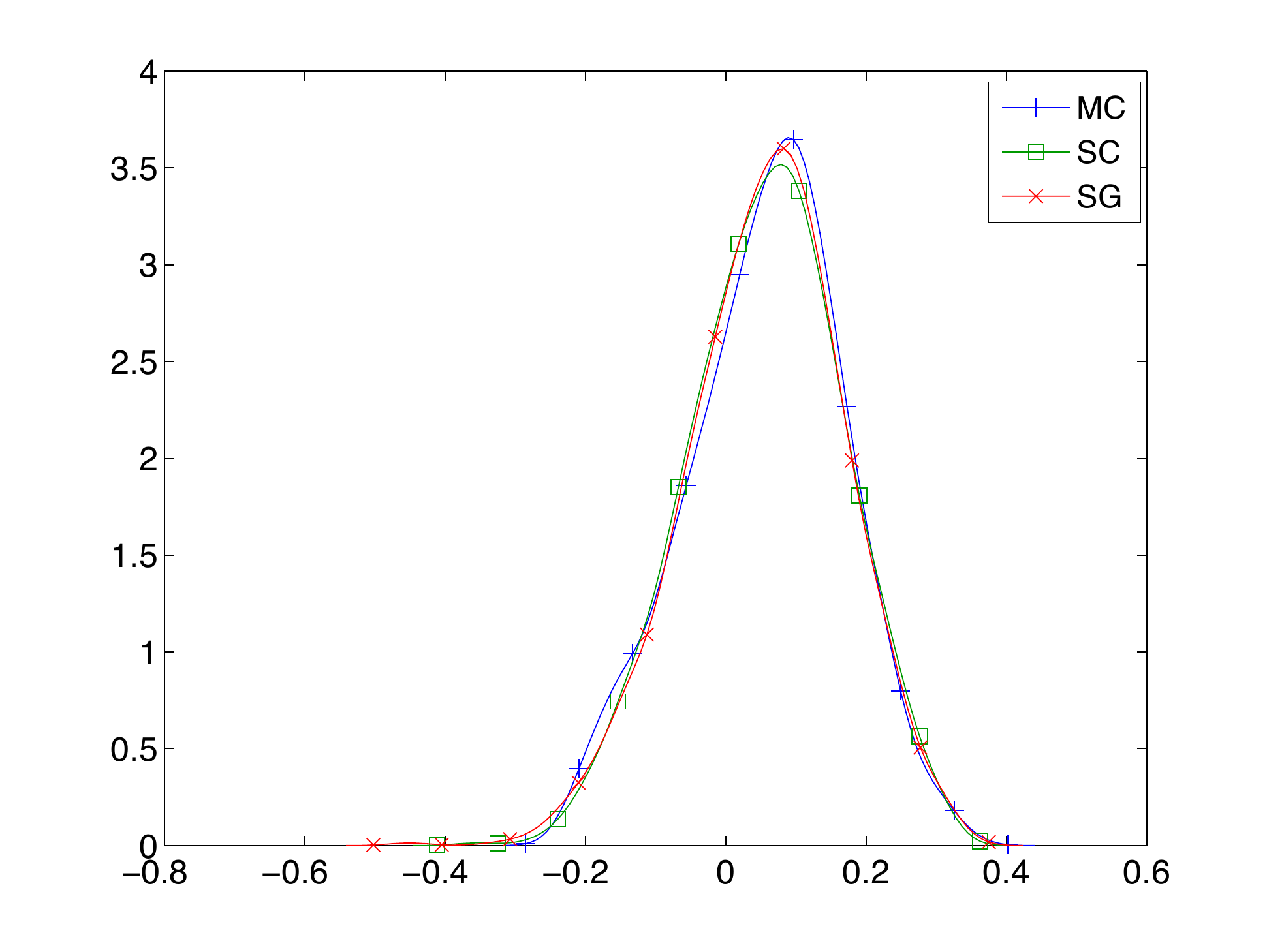}
\end{center}
\caption{Estimated pdfs of the velocities $u_{x}$ (top), $u_{y}$ (middle), and
pressure$ p$ (bottom) with $\operatorname{Re}_{0}=300$ and $CoV=10\%$ (left)
and $30\%$ (right) at the point with coordinates $(3.6436,0)$.}%
\label{fig:sRe300_QoI}%
\end{figure}

As above, we also examined estimated probability density functions for the
velocities and pressures at a specified point in the domain, in this case, at
the point $(3.6436,0)$, taken again from the region in which the solution has
large variance. Figure~\ref{fig:sRe300_QoI} shows these pdf estimates, from
which it can be seen that the three methods for handling uncertainty are in
close agreement for each of the two values of $CoV$.

\begin{figure}[t]
\begin{center}
\includegraphics[width=6.25cm]{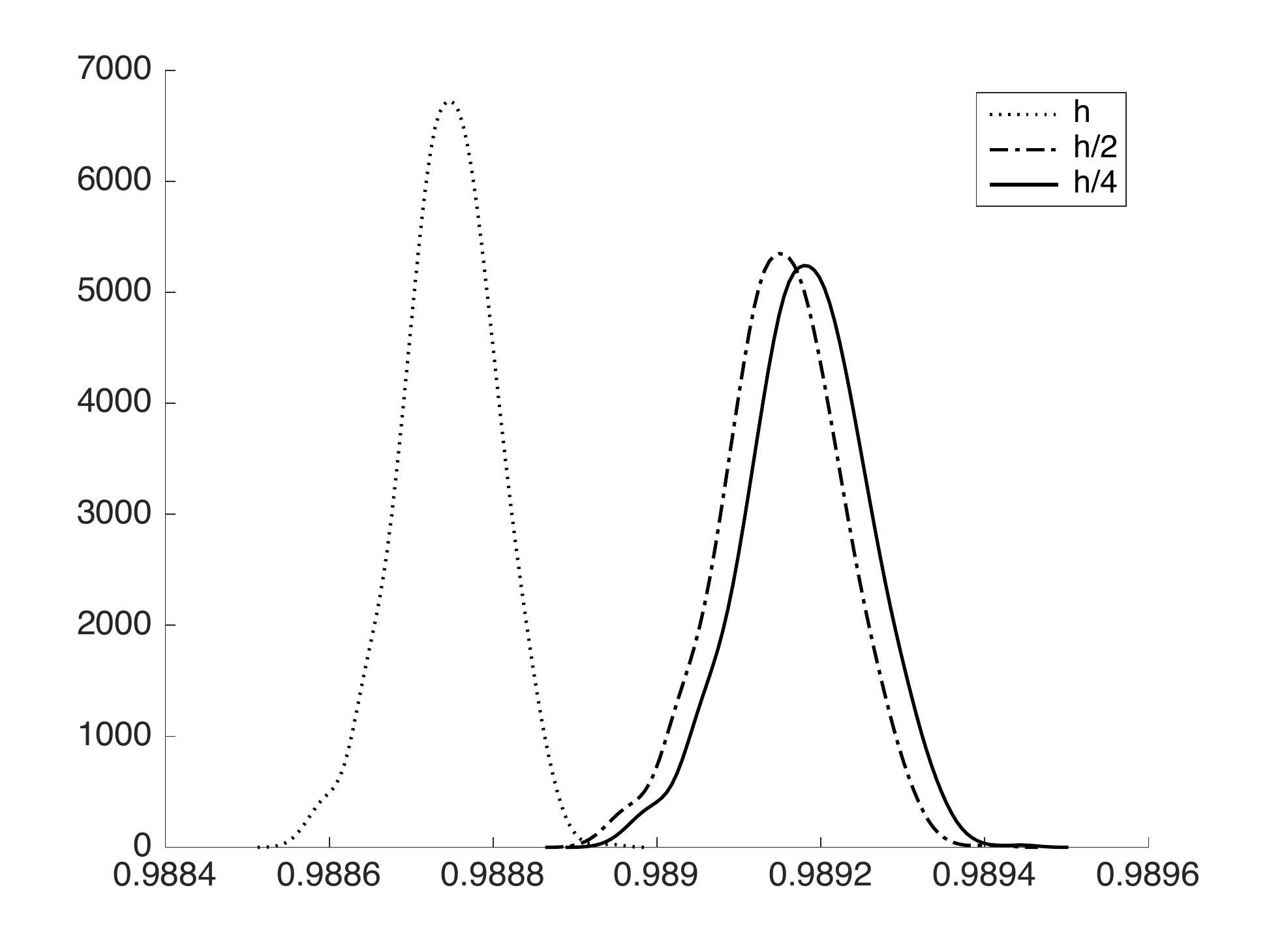}
\includegraphics[width=6.25cm]{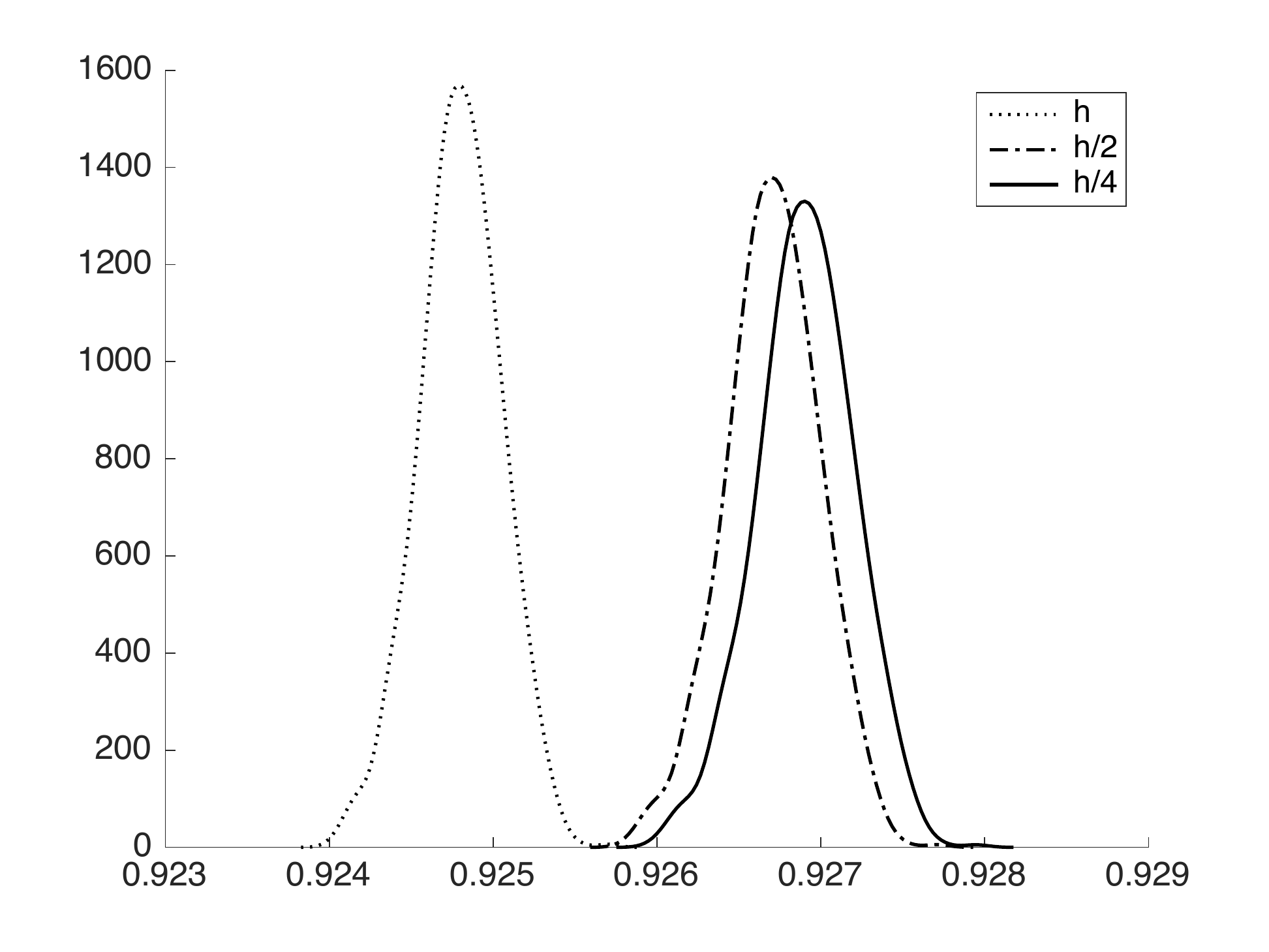}
\vspace{.075in}
\newline\includegraphics[width=6.25cm]{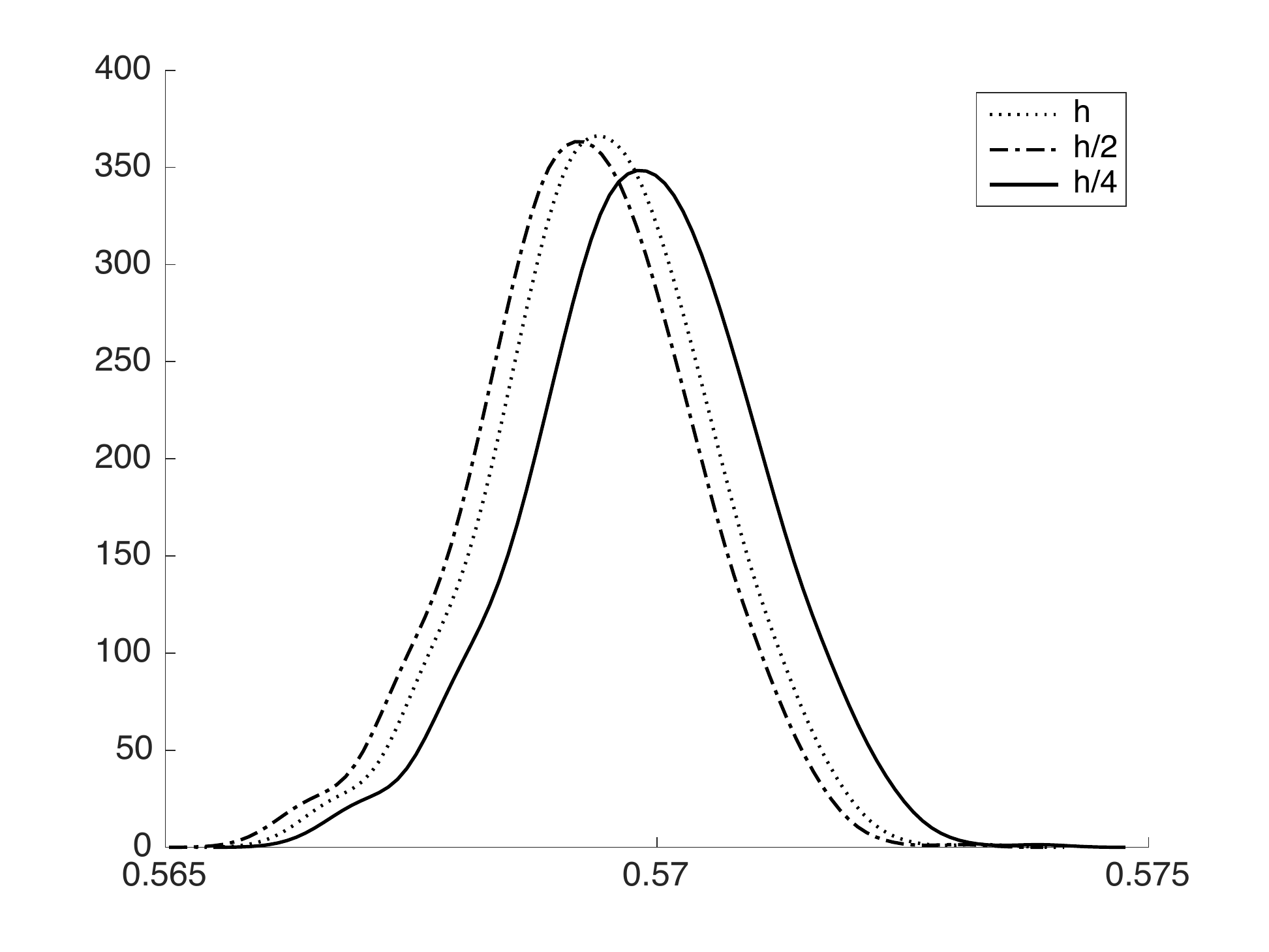}
\includegraphics[width=6.25cm]{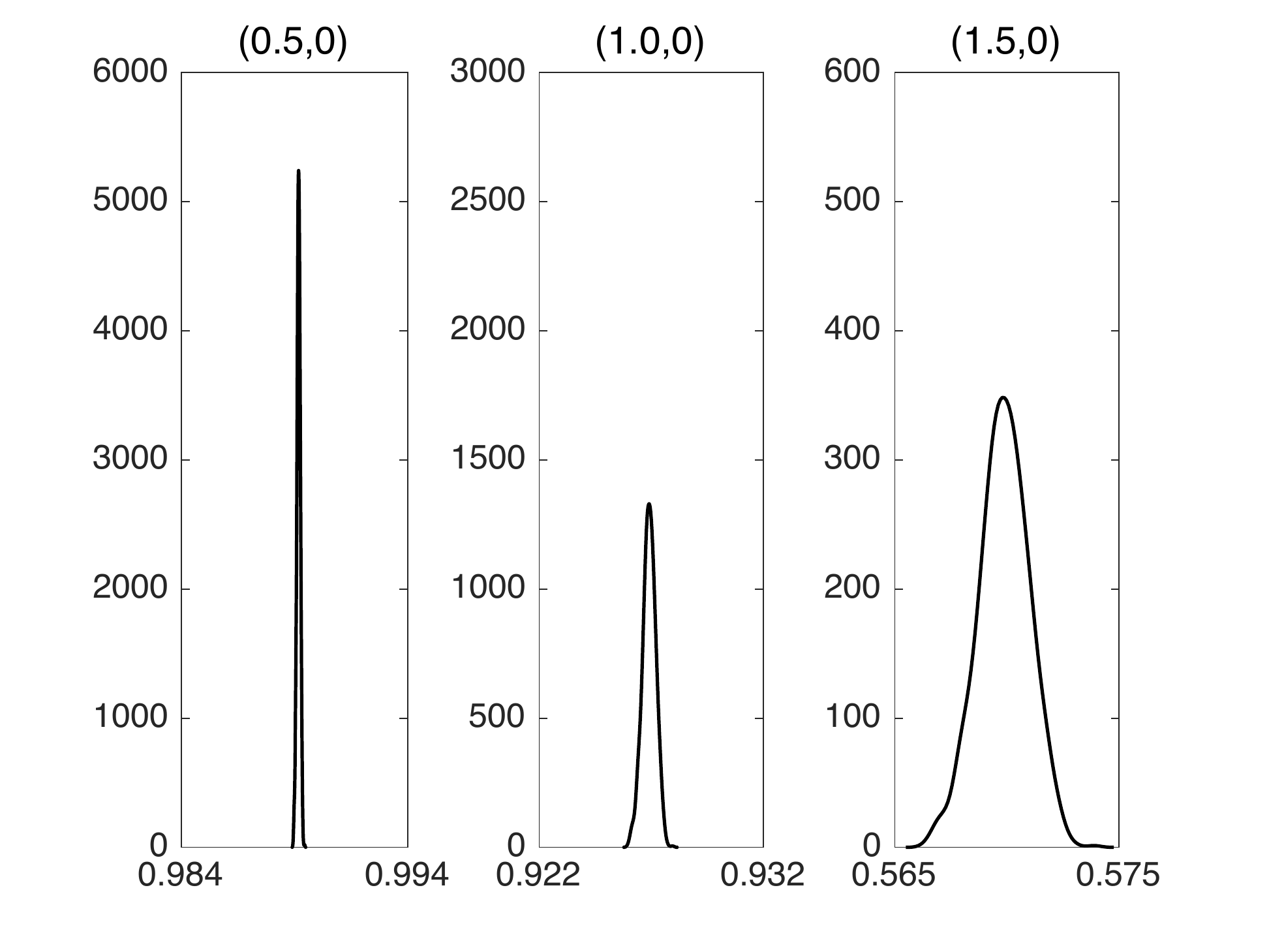}
\end{center}
\caption{Estimated pdfs of the velocity $u_{x}$ at points $(0.5,0)$ (top left)
$(1,0)$ (top right) and $(1.5,0)$ (bottom left), with $\operatorname{Re}_{0}=300$ and $CoV=10\%$ and three meshes.
Bottom right:  detailed depiction of the pdfs for each grid point obtained from the finest mesh.}%
\label{fig:inflow_QoI}%
\end{figure}

\textcolor{black}{Finally, we show in Figure \ref{fig:inflow_QoI} the results of one additional 
experiment with $\operatorname{Re}_0=300$ and $CoV=10\%$, where estimated pdfs of the first velocity component~$u_{x}$ were 
computed at three points near the inflow boundary, $(x,y)=(0.5,0)$, $(1,0)$ and $(1.5,0)$.
These plots show some effects of spatial  accuracy. 
{\color{black}The results for all methods were identical, and so only one representative pdf is shown.}  
The images on the top and bottom left show
results for a uniform mesh of width $h=1/8$
and for two refined meshes in which the horizontal mesh width to the left of the obstacle 
is reduced to $h/2$ and $h/4$.  
The image in the bottom right provides a more detailed view of the fine-grid results;
in this image, the width of the horizontal window is the same ($0.01$) for the three subplots
but the vertical heights are different.
Several trends are evident:
\begin{itemize}
\item 
The pdfs at points nearer the boundary are much narrower. 
\item 
When the spatial accuracy is low, {\color{black} the} methods miss some features of the pdf.
{\color{black}In particular} the 
methods produce a pdf that captures its ``spiky'' nature near the inflow, 
but the location of the spike is not correct.
\end{itemize}
We believe these effects stem from the fact that the inflow boundary is deterministic (with $u_y=1$ 
at $y=0$), and the effects of variability in the viscosity are felt less strongly at points 
near the inflow boundary.
At points further from the deterministic inflow boundary, these effects and differences become 
less dramatic.}


\subsection{Nonlinear solvers}

\label{sec:nonlinear} We briefly comment on the nonlinear solution algorithm
used to generate the results of the previous section. The nonlinear solver was
implemented by modifying the analogue for deterministic systems in IFISS. It
uses a hybrid strategy in which an initial approximation is obtained from
solution of the stochastic Stokes problem~(\ref{eq:Stokes-s}), after which
several steps of Picard iteration (equation~(\ref{eq:linearized-system}) with
$\mathcal{F}_{\ell}$ specified using~(\ref{eq:linearized-operator})
and~(\ref{eq:stoch-Picard})) are used to improve the solution, followed by
Newton iteration ($\mathcal{F}_{\ell}$ from~(\ref{eq:stoch-Newton})). A
convergent iteration stopped when the Euclidian norm of the algebraic
residual~(\ref{eq:R_gPC}) satisfied $\Vert\mathcal{R}^{n}\Vert_{2}\leq
\epsilon\Vert y\Vert_{2}$ where $\epsilon=10^{-8}$ and $\mathbf{{y}}$ is as
in~(\ref{eq:Stokes-s}).

In the experiments described in Section~\ref{sec:obstacle}, we used values of
the Reynolds number, $Re=100$ and $300$, and for each of these, two values of
the coefficient of variation, $CoV=10\%$ and $30\%$. We list here the numbers
of steps leading to convergence of the nonlinear algorithms that were used to
generate the solutions discussed above. Direct solvers were used for the
linear systems; we discuss a preconditioned iterative algorithm in
Section~\ref{sec:preconditioner} below. \begin{tabbing}
\hspace{.75in}
\=$Re=100$, \ \=$CoV=10\%$: \hspace{.2in}
\=$6$ Picard steps \hspace{.05in} \=$1$ Newton step(s) \\
\>$Re=100$, \>$CoV=30\%$: \>6  \>{\color{black}2} \\
\>$Re=300$, \>$CoV=10\%$: \>20 \>1 \\
\>$Re=300$, \>$CoV=30\%$: \>20 \>2
\end{tabbing}
Thus, a larger $CoV$ (larger standard deviation of the random field
determining uncertainty in the process) leads to somewhat larger computational
costs. For $Re=300$, the nonlinear iteration was not robust with $6$ initial
Picard steps (for the stochastic Galerkin method as well as the sampling
methods); $20$ steps was sufficient.

We also explored an \emph{inexact} variant of these methods, in which the
coefficient matrix of~(\ref{eq:linearized-system}) for the Picard iteration
was replaced by the block diagonal matrix $\mathbf{H}_{0}\otimes
\mathcal{F}_{0}^{n}$ obtained from the mean coefficient. For $CoV=10\%$, with
the same number of (now inexact) Picard steps as above ($6$ for $Re=100$ and
$20$ for $Re=300$), this led to just one extra (exact) Newton step for
$Re=100$ and no additional steps for $Re=300$. On the other hand, for
$CoV=30\%$, this inexact method failed to converge.

\section{Preconditioner for the linearized systems}

\label{sec:preconditioner} The solution of the linear systems required during
the course of the nonlinear iteration is a computationally intensive task, and
use of direct solvers may be prohibitive for large problems. In this section,
we present a preconditioning strategy for use with Krylov subspace methods to
solve these systems, and we compare its performance with that of several other
techniques. The new method is a variant of the hierarchical Gauss-Seidel
preconditioner developed in~\cite{Sousedik-2014-THP}.

\subsection{Structure of the matrices and the preconditioner}

\label{sec:matrices} 
\begin{figure}[ptbh]
\centering
\begin{tabular}
[c]{cc}%
\includegraphics[width=5.7cm]{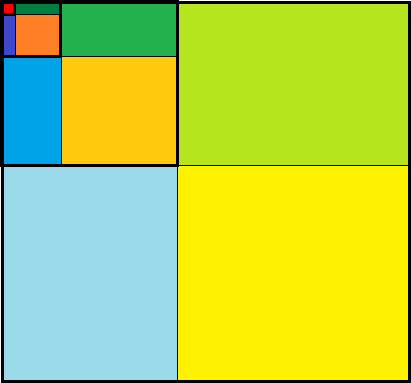} &
\includegraphics[width=5.6cm]{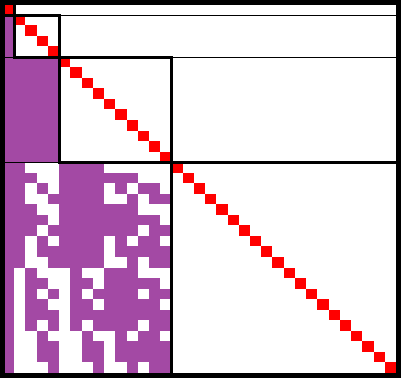}
\end{tabular}
\caption{Hierarchical structure of the stochastic Galerkin
matrix~(\ref{eq:hSG}) (left) and splitting operator~$L+$diag$(D)$\ for the
approximate hierarchical Gauss-Seidel preconditioner (ahGS) with the
truncation of the MATVEC (right).}%
\label{fig:GS}%
\end{figure}

We first recall the structure of the matrices
$\{\mathbf{H}_{\ell}\}$ of~(\ref{eq:stochastic-matrix-forms}). More
comprehensive overviews of these matrices can be found
in~\cite{Ernst-2010-SGM,Matthies-2005-GML}. The matrix structure can be
understood through knowledge of the coefficient matrix $c_{P}\equiv\sum
_{\ell=0}^{M_{\nu}-1}h_{\ell,jk}$ where $j,k=0,\dots,M-1$. The block sparsity
structure depends on the type of coefficient expansion
in~(\ref{eq:viscosity-gPC}). If only linear terms are included, that is
$\psi_{\ell}=\xi_{\ell}$, $\ell=1,\dots,N$, then the coefficients $h_{\ell
,jk}=\mathbb{E}\left[  \xi_{\ell}\psi_{j}\psi_{k}\right]  $ yield a Galerkin
matrix with a block sparse structure. In the more general case, $h_{\ell
,jk}=\mathbb{E}\left[  \psi_{\ell}\psi_{j}\psi_{k}\right]  $ and the
stochastic Galerkin matrix becomes fully block dense. In either case, for
fixed $\ell$ and a set of degree $\mathscr{P}$ polynomial expansions, with
$1\leq\mathscr{P}\leq P$, the corresponding coefficient
matrix~$c_{\mathscr{P}}$ has a hierarchical structure
\[
c_{\mathscr{P}}=\left[
\begin{array}
[c]{cc}%
c_{\mathscr{P}-1} & b_{\mathscr{P}}^{T}\\
b_{\mathscr{P}} & d_{\mathscr{P}}%
\end{array}
\right]  ,\qquad\mathscr{P}=1,\dots,P.
\]
Now, let $\mathcal{A}_{P}$ denote the global stochastic Galerkin matrix
corresponding to either a Stokes problem~(\ref{eq:Stokes-s}) or a linearized
system~(\ref{eq:linearized-system}); we will focus on the latter system in the
discussion below. The matrix~$\mathcal{A}_{P}$ also has a hierarchical
structure
\begin{equation}
\mathcal{A}_{\mathscr{P}}=\left[
\begin{array}
[c]{cc}%
\mathcal{A}_{\mathscr{P}-1} & \mathcal{B}_{\mathscr{P}}\\
\mathcal{C}_{\mathscr{P}} & \mathcal{D}_{\mathscr{P}}%
\end{array}
\right]  ,\qquad\mathscr{P}=1,\dots,P,\label{eq:hSG}%
\end{equation}
where $\mathcal{A}_{0}$\ is the matrix of the mean, derived from~$\nu_{0}$
in~(\ref{eq:viscosity-gPC}). This hierarchical structure is shown in the left
side of Figure~\ref{fig:GS}.

We will write vectors with respect to this hierarchy as
\[
x_{\left(  0:\mathscr{P}\right)  }=\left[
\begin{array}
[c]{c}%
x_{\left(  0\right)  }\\
x_{\left(  1\right)  }\\
\vdots\\
x_{\left(  \mathscr{P}\right)  }%
\end{array}
\right]  ,
\]
where $x_{\left(  q\right)  }$ includes all indices corresponding to
polynomial degree $q$, blocked by spatial ordering determined
by~(\ref{eq:coefficient-order}). With this notation, the global stochastic
Galerkin linear system has the form
\begin{equation}
\mathcal{A}_{P}x_{\left(  0:P\right)  }=f_{\left(  0:P\right)  }.\label{eq:SG}%
\end{equation}

To formulate the preconditioner for~(\ref{eq:SG}), we let
$\widetilde{\mathcal{A}}_{0}$ represent an approximation of~$\mathcal{A}_{0}$
and $\widetilde{\mathcal{D}}_{\mathscr{P}}$ represent an approximation
of~$\mathcal{D}_{\mathscr{P}}$. In particular, let
\begin{equation}
\widetilde{\mathcal{D}}_{\mathscr{P}}=\left[
\begin{array}
[c]{ccc}%
\widetilde{\mathcal{A}}_{0} &  & \\
& \ddots & \\
&  & \widetilde{\mathcal{A}}_{0}%
\end{array}
\right]  ,\label{eq:D-approx}%
\end{equation}
where the number of diagonal blocks is given by $\mathscr{P}$.
We will need the action of the inverse
of~$\widetilde{\mathcal{D}}_{\mathscr{P}}$, or an approximation to it, which
can be obtained using an LU-factorization of~$\widetilde{\mathcal{A}}_{0}$, or
using some preconditioner for~$\widetilde{\mathcal{A}}_{0}$, or using a Krylov
subspace solver. A preconditioner~$\mathcal{P}:w_{\left(
0:\mathscr{P}\right)  }\rightarrow v_{\left(  0:\mathscr{P}\right)  }%
$\ for~(\ref{eq:SG}) is then defined as follows:

\begin{algorithm}
\label{alg:ahGS} [Approximate hierarchical Gauss-Seidel preconditioner (ahGS)]
Solve (or solve approximately)
\begin{equation}
\widetilde{\mathcal{A}}_{0}v_{\left(  0\right)  }=w_{\left(  0\right)  },
\label{eq:ahGS-1}%
\end{equation}
and, for $\mathscr{P}${$=1,\ldots P$}, solve (or solve approximately)
\begin{equation}
\label{eq:ahGS-2}\widetilde{\mathcal{D}}_{\mathscr{P}}v_{\left(
\mathscr{P}\right)  } =\left(  w_{\left(  \mathscr{P}\right)  }-\mathcal{C}%
_{\mathscr{P}} v_{\left(  0:\mathscr{P}-1\right)  }\right)  .
\end{equation}
\end{algorithm}

The cost of preconditioning can be reduced further by truncating the
matrix-vector (MATVEC) operations used for the multiplications by the
submatrices~$\mathcal{C}_{\mathscr{P}}$ in~(\ref{eq:ahGS-2}). The idea is as
follows. The system~(\ref{eq:SG}) can be written as
\begin{equation}
\sum_{j=0}^{M-1}\sum_{\ell=0}^{M_{\nu}-1}h_{\ell,jk}\mathcal{F}_{\ell}%
x_{j}=f_{k},\qquad k=0,\dots,M-1,\label{eq:global-system}%
\end{equation}
and the MATVEC with$~\mathcal{A}_{P}$ is given by
\begin{equation}
v_{j}=\sum_{k=0}^{M-1}\sum_{\ell=0}^{M_{\nu}-1}h_{\ell,jk}\mathcal{F}_{\ell
}u_{k},\label{eq:MAT-VEC}%
\end{equation}
where the indices $j$, $k\in\left\{  0,\dots,M-1\right\}  $ correspond to
nonzero blocks in$~\mathcal{A}_{P}$. The truncated MATVEC is an inexact
evaluation of~(\ref{eq:MAT-VEC}) proposed in~\cite[Algorithm 1]%
{Sousedik-2014-THP}, in which the summation over $\ell=0,\dots.M_{\nu}-1$ is
replaced by summation over a subset $\mathcal{M}_{t}\subseteq\left\{
0,\dots,M_{\nu}-1\right\}  $. Figure~\ref{fig:GS} shows the hierarchical
structure of the matrix and of the ahGS preconditioning operator. Both images
in the figure correspond to the choice $P=3$, so that the hierarchical
preconditioning operation~(\ref{eq:ahGS-1})--(\ref{eq:ahGS-2}) requires four
steps. Because $N=4$, the matrix block size is $M=\left(  \!\!%
\begin{array}
[c]{c}%
N+P\\
P
\end{array}
\!\!\right)  =35$. The block-lower-triangular component of the image on the
right in Figure~\ref{fig:GS} shows the hierarchical structure of the ahGS
preconditioning operator with truncation. For the matrix in the left panel,
$M_{\nu}=\left(  \!\!%
\begin{array}
[c]{c}%
N+2P\\
2P
\end{array}
\!\!\right)  =210$, but the index set$~\mathcal{M}_{t}$ includes terms with
indices at most $M-1$ in the accumulation of sums used for$~\mathcal{C}%
_{\mathcal{P}}$. These two heuristics, approximation by~(\ref{eq:D-approx})
and the truncation of MATVECs, significantly improve the sparsity structure of
the preconditioner, on both the block diagonal (through the first technique)
and the block lower triangle (through the second). In the next section, we
will also consider truncated MATVECs with smaller maximal indices.

\subsection{Numerical experiments}

\label{sec:numerical} In this section, we describe the results of experiments
in which the ahGS preconditioner is used with GMRES to solve the systems
arising from both Picard and Newton iteration.\footnote{The preconditioning
operator may be nonlinear, for example, if the block solves
in (\ref{eq:ahGS-1})--(\ref{eq:ahGS-2}) are performed approximately using
Krylov subspace methods, so that care must be made in the choice of the Krylov
subspace iterative method. For this reason, we used a flexible variant of the
GMRES method~\cite{Saad-1993-FIP}.} 
\textcolor{black}{Unless otherwise specified, in the experiments discussed below, the 
block-diagonal matrices of (\ref{eq:D-approx}), (\ref{eq:ahGS-1}), (\ref{eq:ahGS-2}) are taken 
to be $\widetilde{\mathcal{A}}_{0}=\mathcal{A}_{0}$,
which corresponds to the first summand $\mathbf{H}_{0}\otimes\mathcal{F}_{0}^{n}$ in 
(\ref{eq:linearized-system}) and represents an approximation to the diagonal blocks of 
$\mathcal{D}_{\mathscr{P}}$;
see also Remark \ref{rem:order}. }
We also compare the performance of ahGS preconditioning with several alternatives:
\begin{itemize}
\item Mean-based preconditioner~(MB)~\cite{Pellissetti-2000-ISS,
Powell-2009-BDP}, where the preconditioning operator is the block-diagonal
matrix $\mathbf{H}_{0}\otimes\mathcal{F}_{0}^{n} = I\otimes\mathcal{F}_{0}%
^{n}$ derived from the mean of~$\nu$.

\item The Kronecker-product preconditioner~\cite{Ullmann-2010-KPD}
(denoted~K), given by $\widehat{\mathbf{H}}_{0}\otimes\mathcal{F}_{0}^{n}$,
where $\widehat{\mathbf{H}}_{0}$ is chosen to to minimize a measure of the
difference $\mathcal{A}_{P} - \widehat{\mathbf{H}}_{0}\otimes\mathcal{F}%
_{0}^{n}$.

\item Block Gauss-Seidel preconditioner~(bGS), in which the preconditioning
operation entails applying $(\widetilde{\mathcal{D}}_{\mathscr{P}}+L)^{-1}$
determined using the block lower triangle of~$\mathcal{A}_{\mathcal{P}}$, that
is with the approximation of the diagonal blocks by$~\mathcal{A}_{0}$ but
without MATVEC truncation. This is an expensive operator but enables an
understanding of the impact of various efforts to make the preconditioning
operator more sparse.


\item \textcolor{black}{ahGS(PCD), a modification of ahGS in which the block-diagonal 
matrices of (\ref{eq:ahGS-1})--(\ref{eq:ahGS-2}) are replaced by the 
\emph{pressure convection-diffusion} (PCD) approximation \cite{Elman-2014-FEF} to 
the mean matrix  $\mathcal{F}_{0}$ .}

\item
 \textcolor{black}{ahGS(PCD-it),  in which the block-diagonal solve of (\ref{eq:ahGS-2}) is 
replaced by an approximate solve determined by twenty steps of a PCD-preconditioned 
GMRES iteration.
For this strategy, 
$\widetilde{\mathcal{A}}_{0}$  in (\ref{eq:ahGS-1}) and the blocks of 
$\widetilde{\mathcal{D}}_{\mathscr{P}}$  in (\ref{eq:ahGS-2}) are taken to be the
complete sum of (\ref{eq:linearized-system}), but the approximate solutions to these
systems are obtained using a preconditioned inner iteration.}
\end{itemize}

Four of the strategies, the ahGS, mean-based, Kronecker-product and bGS
preconditioners, require the solution of a set of block-diagonal systems with
the structure of a linearized Navier-Stokes operator of the form given
in~(\ref{eq:Newton}). We used direct methods for these computations. All
results presented are for $\operatorname{Re}=100$; performance for
$\operatorname{Re}=300$ were essentially the same. We believe this is because
of the exact solves performed for the mean operators.

The results for the first step of Picard iteration are in
Tables~\ref{tab:Re_100-Picard-N}--\ref{tab:Re_100-Picard-trunc}. All tests
started with a zero initial iterate and stoppped when the residual
$r^{(k)}=f_{(0:P)}-\mathcal{A}_{P}x_{(0:P)}^{(k)}$ for the $k$'th iterate
satisfied $\Vert r^{(k)}\Vert_{2}\leq10^{-8}\Vert|f_{(0:P)}\Vert_{2}$ in the
Euclidian norm. With other parameters fixed and no truncation of the MATVEC,
Table~\ref{tab:Re_100-Picard-N} shows the dependence of
GMRES iterations on the stochastic dimension~$N$,
Table~\ref{tab:Re_100-Picard-P} shows the dependence on the degree of
polynomial expansion~$P$, and Table~\ref{tab:Re_100-Picard-CoV} shows the
dependence on the coefficient of variation~$CoV$. It can be seen that the
numbers of iterations with the ahGS preconditioner are essentially the same as
with the block Gauss-Seidel~(bGS)\ preconditioner, and they are much smaller
compared to the mean-based~(MB) and the Kronecker product preconditioners. 
When the exact solves with the mean matrix are replaced by the
mean-based modified pressure-convection-diffusion (PCD) preconditioner for the
diagonal block solves, the iterations grow rapidly. 
\textcolor{black}{On the other hand, if the PCD preconditioner is used as part of an 
inner iteration as in ahGS(PCD-it), then good performance is recovered.}
This indicates that a good preconditioner for the mean matrix is an essential component 
of the global preconditioner for the stochastic Galerkin matrix.

Table~\ref{tab:Re_100-Picard-trunc} shows the iteration counts when the MATVEC
operation is truncated in the action of the preconditioner. Truncation
decreases the cost per iteration of the computation, and it can be also seen
that performance can actually be improved. For example, with \thinspace
$\ell_{t}=2$, the number of iterations is the smallest. Moreover there are
only $21$ nonzeros in the lower triangular part of the sum of coefficient
matrices~$\{\mathbf{H}_{\ell}\}$ (each of which has order $10$) used in the
MATVEC with $\ell_{t} \leq2$, compared to $63$ nonzeros when no truncation is
used; there are $203$ nonzeros in the set of $28$ full matrices~$\{\mathbf{H}%
_{\ell}\}$.

\begin{table}[ptbh]
\caption{For Picard step: dependence on stochastic dimension $N$ of GMRES
iteration counts, for various preconditioners, with polynomial degree $P=3$
and coefficient of variation $CoV=30\%$. $M$ is the block size of the
stochastic Galerkin matrix, $M_{\nu}$ the number of terms
in~(\ref{eq:viscosity-gPC}) and $ngdof$ the size of the stochastic Galerkin
matrix.}%
\label{tab:Re_100-Picard-N}
\begin{center}%
\begin{tabular}
[c]{|c|ccc|c|c|c|c|c|c|}\hline
$N$ & $M$ & $M_{\nu}$ & $ngdof$ & MB & K & bGS & ahGS & \textcolor{black}{ahGS(PCD)} & \textcolor{black}{ahGS(PCD-it)} \\\hline
1 & 4 & 7 & 57,120 & 63 & 36 & 30 & 30 & \textcolor{black}{201} & \textcolor{black}{29}\\
2 & 10 & 28 & 142,800 & 102 & 66 & {\color{black}59} & 54 & \textcolor{black}{357} & \textcolor{black}{48}\\
3 & 20 & 84 & 285,600 & 145 & 109 & 88 & 82 & \textcolor{black}{553} & \textcolor{black}{73}\\\hline
\end{tabular}
\end{center}
\end{table}

\begin{table}[ptbh]
\caption{For Picard step: dependence on polynomial degree $P$ of GMRES
iteration counts, for various preconditioners, with stochastic dimension $N=3$
and coefficient of variation $CoV=30\%$. Other headings are as in
Table~\ref{tab:Re_100-Picard-N}.}%
\label{tab:Re_100-Picard-P}
\begin{center}%
\begin{tabular}
[c]{|c|ccc|c|c|c|c|c|c|}\hline
$P$ & $M$ & $M_{\nu}$ & $ngdof$ & MB & K & bGS & ahGS & \textcolor{black}{ahGS(PCD)} & \textcolor{black}{ahGS(PCD-it)}\\\hline
1 & 4 & 10 & 57,120 & 26 & 22 & 12 & 11 & \textcolor{black}{144} & \textcolor{black}{13}\\
2 & 10 & 35 & 142,800 & 63 & 49 & 29 & 26 & \textcolor{black}{300} & \textcolor{black}{29}\\
3 & 20 & 84 & 285,600 & 145 & 109 & 88 & 82 & \textcolor{black}{553} & \textcolor{black}{73}\\\hline
\end{tabular}
\end{center}
\end{table}

\begin{table}[ptbh]
\caption{For Picard step: dependence on coefficient of variation $CoV$ of
GMRES iteration counts, for various preconditioners, with stochastic dimension
$N=2$ and polynomial degree $P=3$. Other headings are as in
Table~\ref{tab:Re_100-Picard-N}.}%
\label{tab:Re_100-Picard-CoV}
\begin{center}%
\begin{tabular}
[c]{|c|c|c|c|c|c|c|}\hline
$CoV (\%)$ & MB & K & bGS & ahGS & \textcolor{black}{ahGS(PCD)} & \textcolor{black}{ahGS(PCD-it)}\\\hline
10 & 16 & 14 & 7 & 7 & \textcolor{black}{186} & \textcolor{black}{10}\\
20 & 36 & 27 & 17 & 16 & \textcolor{black}{259}  & \textcolor{black}{17}\\
30 & 102 & 66 & {\textcolor{black}59} & 54 & \textcolor{black}{357}  & \textcolor{black}{48}\\\hline
\end{tabular}
\end{center}
\end{table}

\begin{table}[ptbh]
\caption{For Picard step: number of GMRES iterations when the preconditioners
use truncated MATVEC. The matrices corresponding to higher degree expansion of
the coefficient than~$\ell_{t}$ are dropped from the action of the
preconditioner. Here $N=2$, $P=3$, $M_{t}$ is the number of terms used in the
inexact (truncated) evaluation of the MATVEC~(\ref{eq:MAT-VEC}), $nnz(c_{\ell
})$ is the number of nonzeros in the lower triangular parts of the coefficient
matrices~(\ref{eq:stochastic-matrix-forms}) after the truncation. Other
headings are as in Table~\ref{tab:Re_100-Picard-N}.}%
\label{tab:Re_100-Picard-trunc}
\begin{center}%
\begin{tabular}
[c]{|cc|c|c|c|c|c|}\hline
\multirow{3}{*}{setup} & $\ell_{t}$ & 0 & 1 & 2 & 3 & 6\\
& $M_{t}$ & 1 & 3 & 6 & 10 & 28\\
& $nnz(c_{\ell})$ & 0 & 12 & 21 & 43 & 63\\\hline
\multicolumn{7}{|c|}{bGS}\\\hline
\multirow{3}{*}{CoV (\%)} & 10 & 16 & 9 & 7 & 7 & 7\\
& 20 & 36 & 19 & 15 & 17 & 17\\
& 30 & 102 & 55 & 43 & 57 & 58\\\hline
\multicolumn{7}{|c|}{ahGS}\\\hline
\multirow{3}{*}{CoV (\%)} & 10 & 16 & 9 & 7 & 7 & 7\\
& 20 & 36 & 19 & 14 & 16 & 16\\
& 30 & 102 & 55 & 35 & 55 & 54\\\hline
\end{tabular}
\end{center}
\end{table}

Results for the first step of Newton iteration, which comes after six steps of
Picard iteration, are summarized in Tables~\ref{tab:Re_100-Newton-N}%
--\ref{tab:Re_100-Newton-trunc}. As above, the first three tables show the
dependence on stochastic dimension $N$ (Table~\ref{tab:Re_100-Newton-N}),
polynomial degree $P$ (Table~\ref{tab:Re_100-Newton-P}), and coefficient of
variation $CoV$ (Table~\ref{tab:Re_100-Newton-CoV}). It can be seen that all
iteration counts are higher compared to corresponding results for the Picard
iteration, but the other trends are very similar. In particular, the
performances of the ahGS and bGS preconditioners are comparable except for the
case when $N=3$, and $P=3$ (the last rows in Tables~\ref{tab:Re_100-Newton-N}
and~\ref{tab:Re_100-Newton-P}). Nevertheless, checking results in
Table~\ref{tab:Re_100-Newton-trunc}, which shows the effect of truncation, it
can be seen that with the truncation of the MATVEC the iteration counts of the
ahGS\ and bGS\ preconditioners can be further improved. Indeed, running one
more experiment for the aforementioned case when $N=3$, and $P=3$, it turns
out that with $\ell_{t}=2$ the number of iterations with the ahGS and bGS
preconditioners are $118$ and $120$, respectively and with $\ell_{t}=3$ they
are $101$ and $104$,\ respectively. That is, the truncation leads to fewer
iterations also in this case, and the performance of ahGS and
bGS\ preconditioners is again comparable.


\begin{table}[ptbh]
\caption{For Newton step: dependence on stochastic dimension $N$ of GMRES
iteration counts, for various preconditioners, with polynomial degree $P=3$
and coefficient of variation $CoV=30\%$. Other headings are as in
Table~\ref{tab:Re_100-Picard-N}.}%
\label{tab:Re_100-Newton-N}
\begin{center}%
\begin{tabular}
[c]{|c|ccc|c|c|c|c|c|c|c|}\hline
$N$ & $M$ & $M_{\nu}$ & $ngdof$ & MB & K & bGS & ahGS & \textcolor{black}{ahGS(PCD)} & \textcolor{black}{ahGS(PCD-it)}\\\hline
1 & 4 & 7 & 57,120 & 73 & 42 & 32 & 32 & \textcolor{black}{309} & \textcolor{black}{40}\\
2 & 10 & 28 & 142,800 & 126 & 80 & 77 & 69 & \textcolor{black}{546} & \textcolor{black}{63}\\
3 & 20 & 84 & 285,600 & 235 & 170 & 128 & 151 & \textcolor{black}{1011} & \textcolor{black}{117}\\\hline
\end{tabular}
\end{center}
\end{table}

\begin{table}[ptbh]
\caption{For Newton step: dependence on polynomial degree $P$ of GMRES
iteration counts, for various preconditioners, with stochastic dimension $N=3$
and coefficient of variation $CoV=30\%$. Other headings are as in
Table~\ref{tab:Re_100-Picard-N}.}%
\label{tab:Re_100-Newton-P}
\begin{center}%
\begin{tabular}
[c]{|c|ccc|c|c|c|c|c|c|}\hline
$P$ & $M$ & $M_{\nu}$ & $ngdof$ & MB & K & bGS & ahGS & \textcolor{black}{ahGS(PCD)} & \textcolor{black}{ahGS(PCD-it)}\\\hline
1 & 4 & 10 & 57,120 & 27 & 23 & 12 & 11 & \textcolor{black}{219} & \textcolor{black}{32}\\
2 & 10 & 35 & 142,800 & 68 & 55 & 33 & 32 & \textcolor{black}{450} & \textcolor{black}{42}\\
3 & 20 & 84 & 285,600 & 235 & 170 & 128 & 151 & \textcolor{black}{1011} & \textcolor{black}{117}\\\hline
\end{tabular}
\end{center}
\end{table}

\begin{table}[ptbh]
\caption{For Newton step: dependence on coefficient of variation $CoV$ of
GMRES iteration counts, for various preconditioners, with stochastic dimension
$N=2$ and polynomial degree $P=3$. Other headings are as in
Table~\ref{tab:Re_100-Picard-N}.}%
\label{tab:Re_100-Newton-CoV}
\begin{center}%
\begin{tabular}
[c]{|c|c|c|c|c|c|c|}\hline
$CoV (\%)$ & MB & K & bGS & ahGS & \textcolor{black}{ahGS(PCD)} & \textcolor{black}{ahGS(PCD-it)}\\\hline
10 & 17 & 15 & 8 & 8 & \textcolor{black}{322} & \textcolor{black}{28}\\
20 & 40 & 30 & 20 & 19 & \textcolor{black}{379} & \textcolor{black}{35}\\
30 & 126 & 80 & 77 & 69 & \textcolor{black}{546} & \textcolor{black}{63}\\\hline
\end{tabular}
\end{center}
\end{table}

\begin{table}[ptbh]
\caption{For Newton step: number of GMRES iterations when the preconditioners
use the truncation of the MATVEC. Headings are as in
Table~\ref{tab:Re_100-Picard-trunc}.}%
\label{tab:Re_100-Newton-trunc}
\begin{center}%
\begin{tabular}
[c]{|cc|c|c|c|c|c|}\hline
\multirow{3}{*}{setup} & $\ell_{t}$ & 0 & 1 & 2 & 3 & 6\\
& $M_{t}$ & 1 & 3 & 6 & 10 & 28\\
& $nnz(c_{\ell})$ & 0 & 12 & 21 & 43 & 63\\\hline
\multicolumn{7}{|c|}{bGS}\\\hline
\multirow{3}{*}{CoV (\%)} & 10 & 17 & 10 & 8 & 8 & 8\\
& 20 & 40 & 22 & 18 & 20 & 20\\
& 30 & 126 & 68 & 57 & 75 & 77\\\hline
\multicolumn{7}{|c|}{ahGS}\\\hline
\multirow{3}{*}{CoV (\%)} & 10 & 17 & 10 & 8 & 8 & 8\\
& 20 & 40 & 22 & 17 & 19 & 19\\
& 30 & 126 & 68 & 45 & 70 & 69\\\hline
\end{tabular}
\end{center}
\end{table}

Finally, we briefly discuss computational costs. For any preconditioner, each
GMRES step entails a matrix-vector product by the coefficient matrix. For
viscosity given by a general probability distribution, this will typically
involve a block-dense matrix, and, ignoring any overhead associated with
increasing the number of GMRES steps, this will be the dominant cost per step.
The mean-based preconditioner requires the action of the inverse of the
block-diagonal matrix $I\otimes\mathcal{F}_{0}^{n}$. This has relatively small
amount of overhead once the factors of $\mathcal{F}_{0}^{n}$ are computed. The
ahGS preconditioner without truncation effectively entails a matrix-vector
product by the block lower-triangular part of the coefficient matrix, so its
overhead is bounded by $50\%$ of the cost of a multiplication by the
coefficient matrix. This is an overestimate because it ignores the
approximation of the block diagonal and the effect of truncation. For example,
consider the case with stochastic dimension $N=2$ and polynomial expansions of
degree $P=3$ for the solution and $2P=6$ for the viscosity; this gives $M=10$
and $M_{\nu}=28$. In Tables~\ref{tab:Re_100-Picard-trunc}
and~\ref{tab:Re_100-Newton-trunc}, $M_{t}$\ indicates how many matrices are
used in the MATVEC operations and $nnz(c_{\ell})$ is the number of nonzeros in
the sum of the lower triangular parts of the coefficient matrices
$\{\mathbf{H}_{\ell}\,|\,\ell=0,\ldots,M_{t}-1\}$. With complete truncation,
$\ell_{t}=0$, and ahGS reduces to the mean-based preconditioner. With no
truncation, $\ell_{t}=6$, and because the number of nonzeros in $\{\mathbf{H}%
_{\ell}\}$ is $203$, the overhead of ahGS is $63/203$, less than $30\%$ of the
cost of multiplication by the coefficient matrix. If truncation is used, in
particular when the iteration count is the lowest ($\ell_{t}=2$), the overhead
is only $21/203\doteq10\%$. Note that with increasing stochastic dimension and
degree of polynomial expansion, the savings will be higher because the ratio
of the sizes of the blocks $\mathcal{C}_{\mathscr{P}}/\mathcal{D}%
_{\mathscr{P}}$ decreases as $P$ increases, see~(\ref{eq:hSG}). Last, the
mean-based preconditioner is embarrassingly parallel; the ahGS preconditioner
requires $P+1$ sequential steps, although each of these steps is also highly
parallelizable. The Kronecker preconditioner is more difficult to assess
because it does not have block-diagonal structure, and we do not discuss it here.

\section{Conclusion}

\label{sec:conclusion} We studied the Navier-Stokes equations with stochastic
viscosity given in terms of polynomial chaos expansion. We formulated the
stochastic Galerkin method and proposed its numerical solution using a
stochastic versions of Picard and Newton iteration, and we also compared its
performance in terms of accuracy with that of stochastic collocation and Monte
Carlo method. Finally, we presented a methodology of Gauss-Seidel hierarchical
preconditioning with approximation using the mean-based diagonal block solves
and a truncation of the MATVEC operations. The advantage of this approach is
that neither the matrix nor the preconditioner need to be formed explicitly,
and the ingredients include only the matrices from the polynomial chaos
expansion and a good preconditioner for the mean-value deterministic problem,
it allows an obvious parallel implementation, and it can be written as a
\textquotedblleft wrapper\textquotedblright\ around existing deterministic code.

\bibliographystyle{plain}
\bibliography{navier-stokes}

\end{document}